\documentclass[10pt,a4paper]{article}
\usepackage{amscd}
\usepackage{latexsym}
\usepackage{amsmath}
\usepackage{amssymb}

\newfont{\bl}{msbm10 scaled \magstep1}
\newfont{\bll}{msbm8}
\begin{document}

\newcommand{\ms}{\medskip\\}
\newcommand{\ep}{\hspace*{\fill}$\Box$}
\newcommand{\R}{\mbox{{\bl R}}}
\newcommand{\N}{\mbox{{\bl N}}}
\newcommand{\NN}{\mbox{{\bll N}}}
\newcommand{\C}{\mbox{{\bl C}}}
\newcommand{\eps}{\varepsilon}

\newcommand{\cA}{{\cal A}}
\newcommand{\cB}{{\cal B}}
\newcommand{\cC}{{\cal C}}
\newcommand{\cE}{{\cal E}}
\newcommand{\cF}{{\cal F}}
\newcommand{\cI}{{\cal I}}
\newcommand{\cL}{{\cal L}}
\newcommand{\cP}{{\cal P}}
\newcommand{\cU}{{\cal U}}

\newcommand{\tX}{\widetilde{X}}
\newcommand{\tT}{\widetilde{T}}
\newcommand{\hsT}{\,^sT}
\newcommand{\hf}{\widehat{f}}
\newcommand{\hg}{\widehat{g}}
\newcommand{\hX}{\widehat{X}}
\newcommand{\hT}{\widehat{T}}
\newcommand{\tr}{\qquad\blacktriangle}

\newcommand{\sa}{\,^\ast a}
\newcommand{\sg}{\,^\ast g}
\newcommand{\ssf}{\,^\ast f}
\newcommand{\sn}{\,^\ast n}
\newcommand{\sA}{\,^\ast A}
\newcommand{\sB}{\,^\ast B}
\newcommand{\sC}{\,^\ast C}
\newcommand{\sF}{\,^\ast F}
\newcommand{\sG}{\,^\ast G}
\newcommand{\sH}{\,^\ast H}
\newcommand{\sK}{\,^\ast K}
\newcommand{\sM}{\,^\ast M}
\newcommand{\sS}{\,^\ast S}
\newcommand{\sT}{\,^\ast T}
\newcommand{\sU}{\,^\ast U}
\newcommand{\sV}{\,^\ast V}
\newcommand{\sW}{\,^\ast W}
\newcommand{\sX}{\,^\ast X}
\newcommand{\sZ}{\,^\ast Z}
\newcommand{\sO}{\,^\ast O}
\newcommand{\sR}{\,^\ast \R}
\newcommand{\sN}{\,^\ast \N}

\def\texts#1{\text{\quad #1 \quad}}

\title{\large \bf NONSTANDARD ANALYSIS IN TOPOLOGY}
\author{*Sergio Salbany\\
{\small \it Mathematics Department, University of South Africa,}\\
{\small \it P.O. Box 392, Pretoria 0003, Republic of South Africa}
\and Todor Todorov\\
{\small \it Mathematics Department, California Polytechnic State University,}\\
{\small \it San Luis Obispo, California 93407, USA}}

\date{}
\maketitle

\begin{abstract} 
We present Nonstandard Analysis by three axioms: the {\em Extension, Transfer and 
Saturation Principles} in the framework of the superstructure of a given infinite 
set. We also present several applications of this axiomatic approach to point-set 
topology. Some of the topological topics such as the Hewitt realcompactification 
and the nonstandard characterization of the sober spaces seem to be new in the 
literature on nonstandard analysis. Others have already close counterparts but 
they are presented here with essential simplifications. 
\end{abstract}
\vfill

*Sergio Salbany has recently passed away after a long illness. He will be remembered as an extraordinary mathematician and a good friend. \\
 
{\large \bf Introduction}\\

Our text consists of three chapters:\\

{\em Chapter  I:} An Axiomatic Approach to Nonstandard Analysis.\\

{\em Chapter  II:} Nonstandard and Standard Compactifications.\\

{\em Chapter  III:} Separation Properties and Monads.
\vspace{3mm}\\

A short description of the contents of these chapters follows:\\

In {\em Chapter I} we present an axiomatic approach to A. Robinson's Nonstandard 
Analysis which is one of the most popular among the researchers applying the 
nonstandard methods as a technique.

Ironically, there are very few expositions based exclusively on these three axioms.
We hope that our text will fill this gap. Although our exposition  is to a large 
extent self-contained, it is not designed for a first introduction to the nonstandard 
theory. Rather, it is written for a reader in mind who has already been through other 
more accessible texts on nonstandard analysis but still lacks the trust and confidence 
needed to apply the nonstandard methods in research. We hope that our text might be 
helpful in this respect. For first reading we recommend the excellent paper by 
Tom Lindstr{\o}m [13], where the nonstandard analysis is presented in terms of sequences, 
equivalence relation and equivalence classes and where, in addition, the reader will 
find a larger list of references on the subject. We should emphasize, however, that 
while the sequential approach, presented in Tom Lindstr{\o}m [13], is, perhaps, 
the best way to start, the axiomatic approach, presented here, is, in our view, 
the best way to apply the nonstandard methods in other fields of mathematics and 
science.

The followers of {\em E. Nelson's Internal Set Theory} [16], who have (finally) decided to 
switch to A. Robinson's framework, are especially warmly welcome. For this special group of 
readers we would like to mention that our attention will be equally directed to both internal 
and external sets; they are both equally important although in somewhat different ways: the 
internal sets are crucial when applying the Transfer and Saturation Principles, while the 
external sets appear in the Extension Principle and in the applications of the nonstandard 
analysis - typically as factor spaces of nonstandard objects.

In {\em Chapter} II and {\em Chapter} III we present some applications of the nonstandard 
methods to point-set topology. However, these topological applications can be also treated 
as {\em exercises} which illustrate and support the theory in Chapter I. We assume a basic 
familiarity with the concepts of point-set topology. We shall use as well the terminology of  
(J.L. Kelley [12]) and (L. Gillman and M. Jerison [3]). For the connection between the 
standard and nonstandard methods in topology we refer to (L. Haddad [5]). We denote by $\N$ 
and $\R$ the sets of the natural and real numbers, respectively, and we also use the notation 
$\N_0 = \{0\}\cup\N$. By $C(X,\R)$ and $C_b(X,\R)$ we shall denote the class of all ``continuous" 
and ``continuous and bounded" functions of the type $f:(X,T)\to(\R,\tau)$, respectively, where  
$(X,T)$ is a topological space and $\tau$  is the usual topology on $\R$.    

Here are more details for these two chapters:

In {\em Chapter} II we describe all Hausdorff compactifications of a given topological space 
$(X,T)$ in the framework of nonstandard analysis. This result is a generalization of an earlier 
work by K.D. Stroyan [21] about the compactifications of completely regular spaces. We 
also give a nonstandard construction of the {\em Hewitt realcompactification}  of a given 
topological space $(X,T)$ which  seems to be new in the literature on nonstandard analysis.

There is numerous works on Hausdorff compactifications of topological spaces in 
nonstandard setting: A. Robinson [17]-[18], W.A.J. Luxemburg [14], M. Machover and J. 
Hirschfeld [15], K.D. Stroyan [21], K.D. Stroyan and W.A.J. Luxemburg [22], H. Gonshor 
[4] and L. Haddad [5] and others. We believe that our description of the Hausdorff 
compactifications, in particular, the Stone-\v Cech compactification is noticeably  simpler 
than those both in the standard and nonstandard literature (mentioned above) mostly due to 
the fact that {\em we manage to avoid the involvement of the weak topology} both on the initial 
space and its compactification. Our technique is based on the concept of the {\em nonstandard 
compactification} $(\sX,\hsT)$ of $(X,T)$, where $\sX$ is the nonstandard extension of $X$  
supplied with the {\em standard topology} $\hsT$, with basic open sets of the form $\sG$, where
$G\in T$. The space $(\sX,\hsT)$ is compact (non Hausdorff), it contains $(X,T)$ densely and every 
continuous function $f$ on $(X,T)$ has a unique continuous extension $\ssf$ on $(\sX,\hsT)$. We 
supply the nonstandard hull $\hX_\Phi = \tX_\Phi / \sim_\Phi$  with the quotient topology $\hT$,
and show that  the space  $(\hX_\Phi,\hT)$ is Hausdorff.  Here,  the set  of  the $\Phi$ - 
{\em finite points} $\tX_\Phi \subseteq \sX$  and the equivalence relation $``\sim_\Phi"$ 
are specified by a family of continuous functions $\Phi$. 
In particular, if $\Phi$ consists of bounded functions only, we have $\tX_\Phi = \sX$ 
and $\hX_\Phi = q[\sX]$, where  $q:\sX\to\hX_\Phi$ is the quotient mapping. Thus, the compactness 
of $\hX_\Phi$  follows simply with the argument that {\em the continuous image of a compact space is 
compact}. In particular, when $\Phi = C_b(X,\R)$ we obtain the Stone-\v Cech compactification 
of $(X,T)$ and by changing $\Phi\subseteq C_b(X,\R)$, we describe in a uniform way all Hausdorff 
compactifications of $(X,T)$. When we choose $\Phi = C(X,\R)$,  we obtain the Hewitt real 
compactification of $(X,T)$. 

We should mention that a similar technique based on the space $(\sX,\hsT)$ has been 
already exploited for studying the {\em compactifications of ordered topological spaces} by the 
authors of this text (S. Salbany, T. Todorov [19]-[20]). 

We should mention as well that the standard topology $\hsT$  is coarser than the {\em discrete}  
$S$-{\em topology} on $\sX$, (known also as $LS$-{\em topology}, where $L$ stands for {\em Luxemburg})
with basic open sets: $^\sigma\cP(X)=\{\sS:S\in\cP(X)\}$, introduced by W.A.J. Luxemburg ([14], p.47 and 
p.55) for a similar purpose. This very property of $\hsT$ allows us to avoid the involvement of 
the weak topology in our construction, thus, to simplify the whole method.

In {\em Chapter} III we study the separation properties of topological spaces such as $T_0,T_1$, 
{\em regularity, normality, complete regularity, compactness and soberness} which are 
characterized in terms of monads. Some of the characterizations have already counterparts in 
the literature on nonstandard analysis (but ours are, as a rule, simpler), while others are 
treated in nonstandard terms for the first time. In particular, it seems that the nonstandard 
characterization of the {\em sober spaces} has no counterparts in the nonstandard literature. We 
also present two new characterizations of the compactness in terms of monads similar to but 
different from $A$. Robinson's famous theorem.\\

\noindent
{\bf{PART I. AN AXIOMATIC APPROACH TO NONSTANDARD ANALYSIS}}

We present Nonstandard Analysis by three axioms: the {\em Extension, Transfer and 
Saturation Principles} in the framework of the superstructure of a given infinite set. We use 
the ultrapower construction only to show the consistency of these axioms. We derive some 
of the basic properties of the nonstandard models needed for the applications presented in the 
next two chapters. Although our exposition  is, to large extend self-contained, it might be 
somewhat difficult for a first introduction to the subject. For first reading we recommend 
Tom Lindstr\o m [13].

\vskip 0,5 truecm
\noindent 
{\bf{1. Preparation of the Standard Theory}}
\vskip 0,3 truecm

In any standard theory the mathematical objects can be classified into two groups: 
{\em abstract points} to which we shall refer as ``standard individuals" (or just ``individuals") and 
sets  (sets of individuals, sets of sets of individuals, sets of sets of sets of 
individuals, etc.). In what follows $S$  denotes the set of the individuals of the standard 
theory under consideration. For example, in Real Analysis we choose $S=\R$, 
in general topology $S=X\cup\R$, where $(X,T)$ is a topological space, in functional analysis
$S =\cal{V}\cup\bf{K}$, where  $\cal{V}$ is a vector space over the scalars $\bf{K}$, etc.

\vskip 0,5 truecm
\noindent 
{\bf{1.1 Definition}} {\em (Superstructure)}: Let $S$ be an infinite set. The superstructure
$V(S)$ on $S$ is the union:
$$
  V(S) = \bigcup_{k\in\N_0} V_k(S),   \leqno{\bf (1.2)}
$$
where $\N_0=\{0\}\cup\N$, $V_k(S)$ are defined inductively by $V_0(S)=S$ and 
$$
  V_{k+1}(S) = V_k(S)\cup\cP(V_k(S)),
$$ 
where $\cP(X)$ denotes the power set of $X$. If $A\in V(S)$, then we define the type $t(A)$
of $A$ by $t(A)= \min\{k\in\N_0:A\in V_k(S)\}$. We shall refer to the 
elements of $V(S)$ as {\em entities} - they are either {\em individuals} if belong to $S$, 
or {\em sets} if belong to $V(S)\setminus S$.

\vskip 0,2 truecm
Notice that 
\begin{eqnarray*}
  S &\;= \;V_0(S)\subset V_1(S)\subset V_2(S)\subset ... ,\\
  S &= \;V_0(S)\in     V_1(S)\in     V_2(S)\in     .... 
\end{eqnarray*}
Hence, it follows that  $V_k(S)\subset V(S)$  and  $V_k(S)\in V(S)$  for all  $k$.

\vskip 0,3 truecm
The most {\em distinguished} property of the superstructure is the {\em transitivity}:

\vskip 0,5 truecm
\noindent 
{\bf{(1.3) Lemma}} {\em (Transitivity)}:  Each $V_k(S)$ is transitive in $V(S)$ in the sense that 
$A\in V_k(S)$ implies either $A\in S$ or $A\subset V_k(S)$. Furthermore, the whole 
superstructure $V(S)$ is transitive (in itself) in the sense that  $A\in V(S)$
implies either $A\in S$, or $A\subset V(S)$.  

\vskip 0,3 truecm
\noindent
{\bf{Proof:}}  $X = V_0(S)$ is obviously transitive. Assume (by induction) that $V_k(S)$
is transitive. Now,  $A\in V_{k+1}(S)$ implies either $A\in V_k(S)$ or $A\subseteq V_k(S)$, 
by the definition of $V_{k+1}(S)$. On the other hand, $A\in V_k(S)$ implies either 
$A\in S$ or $A\subseteq V_k(S)$, by the inductive assumption. 
Hence $V_{k+1}(S)$ is also transitive. The transitivity of the whole $V(S)$ follows immediately: 
$A\in V(S)$ implies $A\in V_k(S)$ for some $k$, thus, we have either $A\in S$ or 
$A\subseteq V_k(S)$, by the transitivity of $V_k(S)$. The latter implies $A\subseteq V(S)$,
since  $V_k(S)\subseteq V(S)$.  $\blacktriangle$

\vskip 0,3 truecm
We observe that the elements of  $S$  are the only elements of  $V(S)$ which are not 
subsets of  $V(S)$. The latter justifies the terminology {\em individuals}  for the 
elements of  $S$. 

The superstructure $V(S)$ consists of all mathematical objects of the  theory: the 
individuals are in $V_0(S)$; the ordered pairs  $\langle x,y\rangle$ in $S\times S$ belongs
to $V_2(S)$ since they can be perceived  as sets of the type $\{\{x\},\{x,y\}\}$;
the functions  $f:S\to S$, and more generally, the relations in $S$ are subsets of 
$V_2(S)$ and hence, belong to $V_3(S)$; if $T$ is a topology on $S$, then $T\subseteq \cP (S)$
and hence $T$ belongs to $V_2(S)$, where $S = X\cup\R$; the algebraic operations in $S$
are perceived as subsets of $S\times S\times S$  and hence also belong to $V(S)$, etc. 

For the study of $V(S)$ we shall use a formal language $\cL(V(S))$ based on bounded 
quantifier formulas : 

\vskip 0,5 truecm\noindent
{\bf{(1.4) Definition}} {\em (The Language $\cL(V(X))$)}:        

\vskip 0,3 truecm 
{\bf{(i)}} The set of the {\em bounded quantifier formulae}  (b.q.f.) $\cL$ consists of the
formulae of the type  $\Phi(x_1,...,x_n)$  that can be made by:  

\vskip 0,3 truecm 
{\bf{a)}} the symbols: 
$=,\in,\neg,\wedge,\vee,\forall,\exists,\Rightarrow,\Leftrightarrow,( ),[ ]$;  

\vskip 0,3 truecm \noindent
and/or

\vskip 0,3 truecm 
{\bf{b)}} countable many variables: $x,y,x_i,A_i,A_j,...,$  etc.;

\vskip 0,3 truecm \noindent
and/or

\vskip 0,3 truecm 
{\bf{c)}} bounded quantifiers of the type $(\forall x\in x_i)$ or $(\exists y\in x_j)$,
$i,j=1,2,...,n$. The variables $x$  and  $y$ are called {\em bounded}  and those which are not
bounded are called {\em free}. 
 
The variables $x_1,...,x_n$ in $\Phi(x_1,...,x_n)$ are exactly the free variables in 
$\Phi(x_1,...,x_n)$.

\vskip 0,3 truecm 
{\bf{(ii)}} Let $S$  be an infinite set and $V(S)$ be superstructure on $S$.
The language $\cL(V(S))$ consists of all statements of the form $\Phi(A_1,...,A_n)$ for some
b. q. f. $\Phi(x_1,...,x_n)\in\cL$ and some $A_1,...,A_n\in V(S)$. The ``points" 
$A_1,...,A_n$ in $\Phi(A_1,...,A_n)$ are called {\em constants} of $\Phi(A_1,...,A_n)$. 

\vskip 0,3 truecm 
The statements in $\cL(V(S))$ can be true or false.

\vskip 0,3 truecm \noindent
{\bf{Warning:}} Formulae including unbounded quantifiers, such as in $(\forall x)(\exists y)$
\newline\noindent
$(x<y)$, {\em are out of} $\cL$!

\vskip 0,3 truecm \noindent
{\bf{(1.5) Example}} {\em (Real Analysis)}: Let $f:\R\to\R$ be a real function in Real Analysis
and let $x_0\in\R$ and $\eps
\in\R_+$. For the set of individuals we choose $S=\R$. Then: 
\begin{eqnarray*}
 &\Phi\,(\eps,x_0,\,f(x_0),\, \R_+,\, \R,\, f,\,<,\,|\,.\,|,\,-)\;= \\
 &= (­\exists\delta\in\R_+)(\forall x\in\R)(\mid x-x_0\mid <\delta\,\Rightarrow\,
  \mid f(x)-f(x_0)\mid < \eps)
\end{eqnarray*}
is a bounded quantifier formula in $\cL(V(\R)$, with {\em constants}:
$\eps,\,x_0,\,f(x_0),\,\R_+,$ $\R,$ $f<,\,\mid\;\mid,``-"$,  
perceived as elements of $V(\R)$ (where $<,\mid\;\mid$ and $``-"$ are the order relation, absolute 
value and subtraction in $\R$, respectively). The above statement might be true or false 
depending on the choice of $\eps,x_0$ and $f$.

For a more detailed exposition of the formal language $\cL(V(S))$ associated with $V(S)$
we refer to (M. Davis [1], Chapter 1) and (T. Lindstr\o m [13], Chapter IV), but we believe 
that the reader can successfully proceed further without a special background in 
mathematical logic.
 
After these preparations of the standard theory we can now involve nonstandard methods.

\vskip 0,5 truecm \noindent
{\bf 2. Axioms of Nonstandard Analsyis}
\vskip 0,3 truecm
 
We present Nonstandard Analysis by means of three axioms (along with the Axiom of 
Choice) known as the {\em Extension, Transfer and Saturation Principles}. The consistences of 
these axioms will be left for the next section.

\vskip 0,3 truecm\noindent
{\bf (2.1) Definiton} {\em (Nonstandard Model)}: Let $S$ be an infinite set 
(of standard individuals for the standard theory under consideration) and $V(S)$ be its 
superstructure. The superstructure $V(\sS)$ of of a given set $\sS$  together with a 
mapping $A\to\sA$ from $V(S)$ into $V(\sS)$ is called a {\em nonstandard model}
of  $S$  if they satisfy the following three axioms:

\vskip0,3 truecm\noindent
{\bf Axiom 1} {\em (Extension Principle)}: $\,^\ast s = s$ for all $s\in S$ or, equivalently, 
$S\subseteq\sS$.

\vskip0,3 truecm\noindent 
{\bf Axiom 2} {\em (Transfer Principle)}: A bounded quantifier formula (b.q.f.)
\newline\noindent
$\Phi(A_1,...,A_n)$ is true in $\cL(V(S))$  iff  its nonstandard counterpart 
$\Phi(\sA_1,...,\sA_n)$ is true in $\cL(V(\sS))$, where $\Phi(\sA_1 ...\sA_n)$ is
obtained from $\Phi(A_1,...,A_n)$ by replacing all constants $A_1,...,A_n$
by their $^\ast$-images $\sA_1,...,\sA_n$,  respectively.

\vskip 0,3 truecm
Axiom 3 will be presented a little later.

\vskip0,3 truecm\noindent
{\bf (2.2) Remark}: Notice that  $\sS$  is the image of $S$ under the mapping $^\ast$. 
Once $\sS$ is found, the superstructure $V(\sS)$ is determined by the formula (1.2),
where $S$ is replaced by $\sS$. The formal language $\cL(V(\sS))$ differs from $\cL(V(S))$
only by its constants: they belong to $V(\sS)$ instead of $V(S)$. Hence the formula 
$\Phi(\sA_1,...,\sA_q)$ is interpreted as a statement about $\sA_1,...,\sA_q$.

\vskip0,3 truecm\noindent 
{\bf (2.3) Example}: Let $S=\R$ and $\Phi$ is the formula in $V(\R)$ given in (1.5), 
then its nonstandard counterpart in  $\cL(V(\sR))$ is given by:  

\begin{eqnarray*}
  &\Phi(\eps,x_0,f(x_0),\sR_+,\sR,\sf, <, \mid\,\mid,-) = \\
  &= (\exists \delta\in\sR_+)(\forall x\in\sR)(\mid x-x_0\mid < \delta\Rightarrow\mid
  \ssf(x)-f(x_0) \mid < \eps), 
\end{eqnarray*}
where the $^\ast$-images $\sR$ and $\sR_+$ (of $\R$ and $\R_+$, respectively)  are 
(by definition) the sets of the nonstandard real numbers and positive nonstandard real numbers, 
respectively, the $^\ast$-image $\ssf$ of $f$ is called (by definition) the 
``nonstandard extension" of f, the asterisks in front of the standard reals are skipped since
$\eps =\,^\ast\eps,$ $x_0 =\,^\ast x_0$ and  $f(x_0)=\ssf(x_0)$, by the Extension 
Principle and, in addition, the asterisks in front of  $^\ast<,\,^\ast\mid\;\mid,\;^\ast-$,
are also skipped, by convention, although these symbols now mean the order relation, 
absolute value and subtraction in $\sR$,  respectively.

\vskip0,5 truecm\noindent 
{\bf (2.4) Definition} {\em (Classification)}:

\vskip0,3 truecm
	{\bf(i)} The entities (individuals or sets) in the range of the $^\ast$-mapping  
are called {\em standard} (although they are actually images of standard objects). 
In other words, $\cA\in V(\sS)$ is  standard  if  $\cA = \sA$  for some  $\in V(S)$. 
If  $A\in V(S)$, then  $\sA$  is called {\em nonstandard extension}  of  $A$. 
Also if  $A\subseteq V(S)$, then the set 
$$
  ^\sigma A = \{\sa : a\in A \}
$$
is called the {\em standard copy} of $A$. In particular,
$$ 
  ^\sigma V(S) = \{\sA : A\in V(S)\}
$$
is the set of all standard entities in  $V(\sS)$. 

\vskip0,3 truecm
	{\bf (ii)} An entity (individual or set) in $V(\sS)$ is called  {\em internal}   
if it is an element of a standard set of  $V(\sS)$. The set of all internal entities is 
denoted by $V_{\mbox{\scriptsize int}}(\sS)$, i.e. 
$$
  V_{\mbox{\scriptsize int}}(\sS) = \{\cA\in V(\sS): \cA\in\sA\quad {\text{for some}}\quad A\in V(S)\}. 
$$
The entities in $V(\sS) - V_{\mbox{\scriptsize int}}(\sS)$ are called {\em external}. 

Notice that the nonstandard individuals in $\sS$  are internal entities. Moreover, if $s\in\sS$, 
then  $s$  is standard (in the sense of the above definition) iff $s\in S$, which justifies 
the terminology {\em standard} introduced above. 

Let  $\kappa$  be an infinite cardinal number. The next (and last) axiom depends on the choice 
of   $\kappa$.

\vskip0,3 truecm\noindent
{\bf Axiom 3} {\em (Saturation Principle: $\kappa$-Saturation)}: $V(\sS)$ is 
$\kappa${\em -saturated}  in the sense that
$$
  \bigcap_{\gamma\in\Gamma}\,\cA_\gamma\not=\emptyset
$$
for any family of internal sets $\{\cA_\gamma\}_{\gamma\in\Gamma}$ in $V(\sS)$  with the finite
intersection property (f.i.p.) and index set  $\Gamma$  with  card  $\Gamma\le\kappa$.

\vskip0,3 truecm\noindent
{\bf (2.8) Definition} {\em (Polysaturation)}: $V(\sS)$ is polysaturated if it is $\kappa$-saturated
for  $\kappa\ge$ card $(V(S))$.

\vskip0,3 truecm\noindent 
{\bf (2.9) Remark} {\em (The Choice of  $\kappa$)}: We should mention that a given standard 
theory  $V(S)$ has actually many nonstandard models $V(\sS)$ although they can be shown to
be isomorphic under some extra set-theoretical assumptions at least in the case when they 
have the same degree of saturation $\kappa$. The choice of $\kappa$, however, is in our hands 
and depends on the standard theory and our specific goals. In particular, if $(X,T)$ 
is a topological space, we apply a $\kappa$-saturated nonstandard model with the set of 
standard individuals $S = X\cup\R$ (a choice $S\supseteq X\cup\R$ is also possible) and a 
degree of saturation $\kappa\ge$ card $\cB$ (or $\kappa\ge$ card $T$), where {\em B}  is a 
base for $T$.

As usual, we can not survive (even in the framework of a superstructure) without the 
axiom of choice:

\vskip 0,3 truecm\noindent 
{\bf Axiom 4} {\em (Axiom of Choice)}:  Let  $I\in V(S)\setminus S$  and  $\{A_i\}_{i\in I}$
be a family of non-empty sets in $V(S)\setminus S$, i.e. $A_i\in V(S)\setminus S$ for all 
$i\in I$. Then there exists a function (of choice) ${\cC} : \to\bigcup_{i\in I}\; A_i$
such that $\cC (i)\in A_i$  for all  $i\in I$. 

\vskip 0,3 truecm\noindent 
{\bf (2.10) Remark}: Although we consider the text presented in this section as an ``up to date" 
version of A. Robinson's Nonstandard Analysis, we should mention that the original A. 
Robinson's theory [17] is based on the ``Enlargement Principle" and the concept for a 
``Countably Comprehensive Model", rather than on the ``Saturation Principle and 
$\kappa$-saturation", as presented here. There exist also other axiomatic formulations of 
nonstandard analysis, e.g. H. J. Keisler [11] axiomatization of $\sR$, the 
``Internal Set Theory", due to E. Nelson[16] and, more recently, C.W. Henson [7] axiomatic 
approach. For a discussion and a general overlook we refer again to Tom Lindstr\o m [13].

\vskip 1 truecm\noindent 
{\bf 3. Existence of Nonstandard Models}

\vskip 0,3 truecm 
The content of this section can be viewed either as a proof of the consistency of Axiom 
1-3 of Nonstandard Analysis, presented in Section 2, or, alternatively, as an independent 
constructive approach to nonstandard analysis.

\vskip 0,3 truecm\noindent   
{\bf (3.1) Theorem} {\em(Consistency)}: For any infinite set  $S$  and any infinite cardinal 
$\kappa$ there exists a $\kappa$-saturated (polysaturated) nonstandard model $V(\sS)$ of $S$.

A sketch of the proof is presented in  A)  and  B)  below. For more detailed exposition 
we refer to T. Lindstr\o m [13]

\vskip 0,3 truecm\noindent 
{\bf A) Existence of $\aleph_0$-Saturated Nonstandard Extensions}: 

\vskip 0,3 truecm
Although Nonstandard Analysis arose historically in close connection with model 
theory and mathematical logic, it is completely possible to construct it in the framework of 
Standard Analysis, i.e. assuming the axioms of Standard Analysis only (along with the Axiom 
of Choice). The method is known as ``ultrapower construction" or ``constructive nonstandard 
analysis". This part of our exposition can be viewed either as a proof of the consistence 
theorem above in the particular case $\kappa = \aleph_0$, where $\aleph_0 =$ card$\;\N$, 
or as an independent ``sequential approach" to Nonstandard Analysis:

\vskip 0,3 truecm
{\bf (i)} Let $p:\cP(\N)\to\{0,1\}$ be a finitely additive measure such that
\newline\noindent 
$p(A)= 0$ for all finite $A\subset \N$ and  $p(\N)=1$. To see that there exist measures 
with these properties, take a free ultrafilter $\cU\subset\cP(\N)$ on $\N$ (here the Axiom
of Choice is involved) and define  $p(A)= 0$ for  $A\notin \cU$ and $p(A)=1$ for $A\in\cU$.
We shall keep $p$ fixed in what follows. 

\vskip 0,3 truecm

{\bf (ii)} Let $S^{\Bbb N}$ be the set of all sequences in $S$. Define an equivalence relation 
$\sim$ in $S^{\Bbb N}$  by: $\{a_n\}\sim\{b_n\}$  if  $a_n = b_n$  a. e., where ``a. e." 
stands for ``almost everywhere", i.e. if $p(\{n:a_n = b_n\})=1$. Then the factor space
$\sS = S^{\Bbb N} /\sim$ defines a set of nonstandard individuals. 
(Notice that $\sS$ depends on the choice of the measure $p$.) We shall denote by  
$\langle a_n\rangle$ the equivalence class determined by the sequence $\{a_n\}$. 
The inclusion $S\subset\sS$ is defined by $s\to\langle s, s,...,\rangle$. We can determine 
now the superstructure $V(\sS)$ by (1.1), where $S$ is replaced by $\sS$, and the latter is 
treated as a set of individuals (although it is, actually, a set of sets of sequences). 

\vskip 0,3 truecm
{\bf (iii)} Let $V(S)^{\Bbb N}$  be the set of all sequences in $V(S)$ (i.e. sequences of 
points in $S$, sequences of subsets of $S$, sequences of functions, sequences of 
``mixture of points and functions", ..., sequences of ``everything"). A sequence $\{A_n\}$ 
in $V(S)^{\Bbb N}$ is called ``tame" if there exists $m$ in $\N_0$ such that 
$A_n\in V_m(S)$  for all $n\in\N$ (or, equivalently, for almost all $n$ in $\N$). If 
$\{A_n\}$  is a tame sequence in $V(S)^{\Bbb N}$, then its type $t(\{A_n\})$  is defined 
as the (unique) $k\in \N_0$  such that  $t(A_n) = k$  a.e., where $t(A_n)$  is the type of 
$A_n$  in $V(S)$ defined in $1^\circ$.  To any tame sequence $\{A_n\}$ in $V(S)^{\Bbb N}$ we 
associate an element  $\langle A_n\rangle$  in  $V(\sS)$  by induction on the type of 
$\{A_n\}$: If  $t(\{A_n\})= 0$, then  $\langle A_n\rangle$  is the element in $\sS$, 
defined in (ii). If  $\langle B_n\rangle$  is already defined for all tame sequences 
$\{B_n\}$ in  $V(S)^{\Bbb N}$  with  $t(\{B_n\}) < k$  and  $t(\{A_n\}) = k$, then
$$
  \langle A_n\rangle = \left\{\langle B_n\rangle : \{B_n\}\in V(S)^\Bbb N; t(\{B_n\})<k;  
  \; B_n\in A_n \ \mbox{a.e.}\;\right\}.
$$
The element $\cA\in V(\sS)$ is called ``internal" if it is of the type $\cA= \langle A_n\rangle$ 
for some tame sequence $\{A_n\}$  in  $V(S)^\Bbb N$. The elements of  $V(\sS)$ of the type 
$\sA = \langle A,A,...\rangle$  for some  $A\in V(S)$, are called ``standard". Now we define 
the $^\ast$- mapping  $A\to \sA$  from  $V(S)$  into  $V(\sS)$  and the construction of the 
nonstandard model is complete. We shall leave to the reader to check that this model satisfies 
Axiom 1, Axiom 2 and Axiom 3 for $\kappa = \aleph_0$  treated now as theorems 
(Tom Lindstr\o m [13]). 

\vskip 0,5 truecm\noindent
{\bf B) Existence of  $\kappa$-Saturated Nonstandard Extensions}

\vskip 0,3 truecm

In the case of a general cardinal $\kappa$, a similar construction and proofs to the presented in 
A) can be carried out replacing  $\N$  with an index set  $I$  of cardinality $\kappa$, 
and a $\{0,1\}$ - valued measure on  $\cP(I)$  which is  $\kappa$ - good in the sense 
explained in $T$. Lindstr{\o}m [13], where  $\kappa$  is the successor of  $\kappa$. 
Notice that every measure on  $\cP(\N)$  given by a nonprinciple ultrafilter on  $\N$  is 
$\aleph_1$-good, so this condition `` to be $\kappa$ - good" is not needed explicitly in the case 
$\kappa = \aleph_0$.\\

{\bf 4. Some Basic Properties of  the Nonstandard Models}

\vskip 0,3 truecm

Let (as before)  $S$  be an infinte set and  $V(\sS)$  be a nonstandard model of $S$ in the 
sense of Section 2. We shall study some very basic properties of  $V(\sS)$ with focus on the 
standard and internal entities (Definition 2.4).

\vskip 0,5 truecm\noindent
{\bf (4.1) Lemma} {\em (Internal Entities and Transitivity)}: 

\vskip 0,3 truecm
{\bf (i)} $V_{\mbox{\scriptsize int}}(\sS)$ is a countable union  of $\sV_k(S)$:
$$
  V_{\mbox{\scriptsize int}}(\sS) =  \bigcup_{k\in\NN_0} \sV_k(S)\,.
$$

\vskip 0,2 truecm         
{\bf (ii)} Each $\sV_k(S)$ is transitive in  $V(\sS)$  in the sense that  $A\in\sV_k(S)$
implies either  $A\in S$  or  $A\subset \sV_k(S)$. Furthermore, the whole set  $V_{\mbox{\scriptsize int}}(\sS)$
is {\em transitive} in  $V(\sS)$  in the sense that  $A\in V_{\mbox{\scriptsize int}}(\sS)$  implies either 
$A\in \sS$, or $A\subseteq V_{\mbox{\scriptsize int}}(\sS)$.

\vskip 0,3 truecm\noindent
{\bf Proof}: (i) Assume that  $\cA\in V_{\mbox{\scriptsize int}}(\sS)$, i.e. $\cA\in \sA$  for some $\cA\in V(S)$.
That is  $A\in V_k(S)$  for some $k\in\N_0$, which implies  $A\subseteq V_k(S)$, by the 
transitivity of  $V_k(S)$. It follows  $\sA\subseteq V_k(S)$, by {\em Transfer Principle}, 
hence  $\cA\in \sV_k(S)$. Conversely,  $\cA\in\sV_k(S)$  for some $k$ implies $\cA\in V_{\mbox{\scriptsize int}}(\sS)$,
by the definition of $V_{\mbox{\scriptsize int}}(\sS)$, since  $V_k(S)\in V(S)$.  

\vskip 0,5 truecm
(ii) To show the transitivity, observe that 
$$
  (\forall A\in V_k(S))[A\in S\vee(A\subseteq V_k(S)]
$$
is true in  $\cL(V(S))$, by the transitivity of  $V_k(S)$  (Lemma 1.3). Hence
$$
  (\forall A\in\sV_k(S))[A\in\sS\vee(A\subseteq\sV_k(S)]
$$
is true in  $\cL(V(\sS))$, as required, by Transfer Principle. $\tr$

\vskip 0,5 truecm
{\bf (4.2) Theorem} {\em (Boolean Properties)}: The extension mapping  $A\to\sA$  from
$V(S)$  into  $V(\sS)$  is injective and its restriction  
$$
  ^\ast:V(S)\setminus S\to V(\sS)\setminus \sS 
$$
preserves the Boolean operations,  i.e. if  $A, B \in V(S)\setminus S$, then
\begin{eqnarray*}
  &^\ast(A\cup B)   &=  \sA\cup \sB \\
  &^\ast(A\cap B)   &=  \sA\cap \sB \\
  &^\ast(A\setminus B) &=  \sA\setminus \sB. 
\end{eqnarray*}

\vskip 0,2 truecm\noindent
{\bf Proof}: To show the extension mapping is injective, assume that  $\sA = \sB$  
for some  $A,B\in V(S)$.  That means that the formula  $\Phi(\sA,\sB)=[\sA=\sB]$
is true in  $\cL(V(\sS))$, by Transfer Principle. Hence,  $\Phi(A,B)=[A=B]$  is true in 
$L(V(S))$, by Transfer Principle, i.e.  $A=B$, as required. For the preservation of the 
Boolean operations, suppose, say, that  $A\cup B=C$   for some  $A, B, C\in V(S)\setminus S$. 
We have to show that  $\sA\cup \sB=\sC$. We have  $A, B, C\in V_k(S)$  for some  $k\in\N$ 
(by the definition of  $V(S))$. On the other hand, we have  $A,B,C\subset V_k(S)$, by the 
transitivity of  $V_k(S)$. Now, the equality $A\cup B=C$   can be formalized by the formula:
\begin{eqnarray*}
   \Phi(A,B,C)\;= &[\,(\forall x\in V_k(S))((x\in A)\,\vee\,(x\in B))\Rightarrow (x\in C)\,]\\ 
	          &\wedge\,[\,(\forall z\in V_k(S))((z\in C)\Rightarrow ((z\in A)\vee (z\in B))\,]
\end{eqnarray*}
which is true in $\cL(V(S))$. It follows that its nonstandard version:
\begin{eqnarray*}
  (\sA,\sB,\sC)\;= &[\,(\forall x\in\sV_k(S))((x\in\sA)\,\vee\,(x\in\sB))\,\Rightarrow\,
                                                                           (x\in\sC)\,]\\ 
		   &\wedge\,[\,(\forall z\in\sV_k(S))((z\in\sC)\,\Rightarrow\,
                                                                 ((z\in\sA)\,\vee\,(z\in\sB))\,]
\end{eqnarray*}
is true in $\cL(V(\sS))$, by the Transfer principle. Hence, $\sA\cup\sB=\sC$, as required. 
The preservation of the rest of the Boolean properties is checked similarly.
$\tr$

\vskip 0,5 truecm\noindent
{\bf (4.3) Definition} {\em (Canonical Imbedding)}: If  $A\in V(S)\setminus S$, 
then the injective imbedding 
$$
  A\subseteq \sA, 
$$
defined by  $a\to\sa$,  is called {\em canonical}.

Notice that  $a\in A$  iff  $\sa\in \sA$, by {\em Transfer Principle}, hence this mapping is well 
defined. In addition, it is injective, by the above theorem, which justifies the above definition. 
This imbedding justifies also the terminology {\em nonstandard extension} for  $\sA$.
Notice that the range of this mapping is exactly  $^\sigma A$  (by the definition of $^\sigma A$).
Later in this section we shall show that  $\sA$  is a proper extension of  $^\sigma A$, 
hence it is a proper extension of  $A$ (in the sense of the above imbedding), whenever  $A$
is an infinite set. 

\vskip 0,5 truecm\noindent
{\bf (4.4) Lemma} {\em (Definable Sets)}: Let  $\Phi(x,x_1,x_2,...,x_n)\in\cL$  be a b.q.f. and  
$B,A_1,...,A_n\in V(S)$. Then:
\begin{align*} 
   ^\ast\{x\in B:\,&\Phi(x,A_1,...,A_n)\quad{\text{is true in}}\quad\cL(V(S))\}\;= \\ 
		    &=\{x\in\sB : \Phi(x,\sA_1,...,\sA_n)\quad{\text{is true in}}\quad\cL(V(\sS))\}. 
\end{align*}

\vskip 0,2 truecm\noindent
{\bf Proof}: Denote  
$$
  A = \{x\in B: \Phi(x,A_1,...,A_n)\quad{\text{is true in}}\quad \cL(V(S))\} 
$$
and let  $\sA$  be the nonstandard extension of  $A$. We have to show that  
$$
  \{x\in\sB : \Phi(x,\sA_1,...,\sA_n)\quad{\text{is true in}}\quad \cL(V(\sS))\} = \sA.
$$
Suppose (for contradiction) that 
\begin{align*}
  (\exists x\in\sA)\;&(\Phi(x,\sA_1,...,\sA_n)\quad{\text{is false in}}\quad \cL(V(\sS))\;\vee\\
                     &(\exists x\in\sB\setminus\sA)(\Phi(x,\sA_1,...,\sA_n)
                      \quad{\text{is true in}}\quad \cL(V_{\mbox{\scriptsize int}}(\sS)).
\end{align*}
We have  $\sB\setminus\sA =\,^\ast(B\setminus A)$, by the Boolean properties. 
As a result, the above formula becomes
\begin{align*}
  (\exists x\in\sA)\;&(\Phi(x,\sA_1,...,\sA_n)\quad{\text{is false in}}\quad 
                            \cL(V_{\mbox{\scriptsize int}}(\sS))\;\vee\\ 
	             &(\exists x\in\,^\ast(B\setminus A))(\Phi(x,\sA_1,...,\sA_n)\quad
                      {\text{is true in}}\quad \cL(V_{\mbox{\scriptsize int}}(\sS)).
\end{align*}
This statement is equivalent to
\begin{align*}
   (\exists x\in A)\;&(\Phi(x,A_1,...,A_n)\quad{\text{is false in}}\quad \cL(V(S))\;\vee\\  
		     &(\exists x\in B\setminus A)(\Phi(x,A_1,...,A_n)\quad{\text{is true in}}\quad
                       \cL(V(S)),
\end{align*}
by the {\em Transfer Principle}. The latter contradicting the choice of $A.\qquad \blacktriangle$

\vskip 0,3 truecm\noindent
{\bf (4.5) Examples} {\em (Standard Intervals in $\sR$)}: Let  $a, b\in\R$, $a<b$. 
Let $S=\R$ and $V(\sR)$ be a nonstandard model of $\R$. We have
\begin{eqnarray*}
		&^\ast(a,b) = \{x\in\sR : a < x < b \},\\
                &^\ast[a,b] = \{x\in\sR : a\le x\le b\},\\
                &^\ast[a,b) = \{x\in\sR : a\le x < b \},
\end{eqnarray*}
by the above lemma (applied for  $\Phi(x,a,b) = \{a < x < b\}$  for the first case and similar 
for the others). Notice that the above subsets of $\sR$ are intervals - open, 
closed and semi-open, respectively - in  the order relation in $\sR$.

\vskip 0,3 truecm\noindent
{\bf (4.6) Theorem} {\em (Finite Sets)}: 

\vskip 0,3 truecm
{\bf (i)} If $A\in V(S)\setminus S$  is a finite set, then  $\sA =\,^\sigma A$.  In particular, 
$$
  ^\ast\{a\} = \{\sa\}
$$
for any  $a\in V(S)$.

\vskip 0,3 truecm
{\bf (ii)} If $A\subseteq S$  is a finite set, then  $\sA = A$.

\vskip 0,3 truecm\noindent
{\bf Proof}: (i) We start with the case of a singlet. There exists  $k\in\N$ such that 
$a\in V_k(S)$  which is equivalent to  $\sa\in\sV_k(S)$ (for the same $k$), 
by {\em Transfer Principle}. We observe now that $\{a\}$  can be described as a definable set:
$$
  \{a\} = \left\{x\in V_k(S) : x = a \quad {\text{in}} \quad \cL(V(S))\right\},
$$
which implies
$$
  ^\ast\{a\} = \left\{x\in \sV_k(S) : x = \sa \quad {\text{in}} \quad \cL(V(\sS))\right\},
$$
by the above lemma, applied for  $\Phi(x,a) = [x = a]$. The right hand side of the above formula 
is (obviously) $\{\sa\}$, thus, $\{\sa\} =\,^\ast\{a\}$, as required. 
In the case of an arbitrary finite set $A$, the result follows from the Boolean properties 
of the extension mapping:
$$
  \sA =\,^\ast\bigl(\,\bigcup_{a\in A}\;\{a\}\bigr)\;=\;  
  \bigcup_{a\in A}\;^\ast\{a\}\;=\; \bigcup_{a\in A}\;\{\sa\} =\,^\sigma A,
$$
as required.

\vskip 0,3 truecm 
(ii) follows from (i) since  $^\sigma A = A$, by {\em Extension Principle}.$\tr$

\vskip 0,3 truecm\noindent
{\bf (4.7) Theorem} {\em (Nonstandard Extensions)}: Let  $A\in V(S) \setminus S$
be a set in the superstructure, $^\sigma A$  be its standard image and  $\sA$ be its 
nonstandard extension. Then:

\vskip 0,3 truecm
{\bf (i)}   $\quad  \sA \cap\,^\sigma V(S) =\,^\sigma A$. \newline

{\bf (ii)}  $\quad ^\sigma A \subseteq \sA$.  \newline

{\bf (iii)} $\quad ^\sigma A = \sA \quad$  iff  $\quad A \quad$  is a finite set.

\vskip 0,5 truecm\noindent
{\bf Proof}: (i) $(\subseteq)$ Suppose  $\alpha\in\sA \cap\,^\sigma V(S)$. On one hand,  
$\alpha\in\,^\sigma V(S)$ means  $\alpha = \sa$, for some  $a\in V(S)$. On the other hand,
$\sa\in\sA$  is equivalent to  $a\in A$, by {\em Transfer Principle}. 

\vskip 0,3 truecm
$(\supseteq)$  Suppose now that  $\alpha\in\,^\sigma A$,  i.e.  $\alpha = \sa$  for some 
$a\in A$. On one hand, $\alpha\in\,^\sigma A$  implies  $\alpha\in\,^\sigma V(S)$, since  
$^\sigma A\subset\,^\sigma V(S)$. On the other hand, $a\in A$ is equivalent to  $\sa\in\sA$,
by {\em Transfer Principle},  thus,  $\alpha\in\sA\cap\,^\sigma V(S)$, as required.

\vskip 0,3 truecm
(ii) follows directly from (i). 

\vskip 0,3 truecm
(iii) $(\Leftarrow)$  was shown in Theorem 4.6. $(\Rightarrow)$ Assume that  $A$  is an 
infinite set. We have  $\sA \cap\,^\sigma V(S)=\,^\sigma A$  and  $^\sigma A\subseteq\sA$, 
by (i) and (ii) (just proved). Consider first the case $A = \N$  which implies 
$\sN\cap V(S)=\,^\sigma\N$ and  $^\sigma \N\subseteq\sN$. We want to show that 
$\sN\setminus\,^\sigma\N\not= \emptyset$.  Observe that if $n\in\N$, then the set 
$\sN\setminus\{\sn\}$ is internal (actually, standard), since
$$
  \sN\setminus\{\sn\} = \sN\setminus\,^\ast\{n\} =\,
  ^\ast(\N\setminus\{n\})\in\,^\sigma V(S)\subset V_{\mbox{\scriptsize int}}(\sS).
$$
The family of internal sets  $\{\sN\setminus\{\sn\}\}_{n\in\NN}$  has (obviously) the finite 
intersection property, since $\sN$ is an infinite set. It follows, by {\em Saturation Principle},
that its intersection is not empty,  i.e. $\sN\setminus\,^\sigma\N\not=\emptyset$ (as promised). 
We return to the general case of an infinte set $A$. Without loss of generality we might assume 
that  $\N\subset A$, $\N\not= A$.  The latter implies both  $^\sigma\N\subset\,^\sigma A$, 
$^\sigma\N\not=\,^\sigma A$  and  $\sN\subset\sA$, $\sN\not=\sA$. Suppose (for contradiction) 
that  $^\sigma A = \sA$. By intersecting both sides by  $\sN$, we get  $\sN =\,^\sigma A\cap\sN$.
For the right hand side we have  $^\sigma A\cap\sN\subseteq\,^\sigma V(S)\subset\sN =\,^\sigma\N$,
by (i), hence,  $^\sigma\N = \sN$, a contradiction.$\qquad \blacktriangle$

\vskip 0,3 truecm\noindent
{\bf (4.8) Corollary} {\em (Standard vs. Nonstandard Individuals)}: Let $A\subset S$.  Then:

\vskip 0,2 truecm
{\bf(i)} $\quad\sA\cap S  =  A$.

\vskip 0,2 truecm 
{\bf (ii)} $\quad\; A\subseteq\sA$.

\vskip 0,2 truecm
{\bf (iii)} $\quad A = \sA\;$  iff  $A$  is a finite set. In particular, $S$  and $V(S)$  
            are proper subsets of  $\sS$  and  $V(\sS)$, respectively.

\vskip 0,3 truecm\noindent
{\bf Proof}: We have  $A =\,^\sigma A$  since  $a = \sa$  for all  $a\in A$, by the 
{\em Extension Principle}.  Hence the result follows directly from the previous theorem. 
In particular for  $A = S$, we have  $S\subset\sS$, $S\not=\sS$,  since  $S$  is an infinite set. 
The latter implies  $V(S)\subset V(\sS)$,  $V(S)\not= V(\sS)$. $\tr$

\vskip 0,5 truecm\noindent
{\bf(4.9) Examples} {\em (Real Numbers)}: Let us consider the important particular case  $S = \R$.
The nonstandard individuals are the nonstandard real numbers $\R$. It follows that $\sR$ 
is a proper extension of $\R$, $\R\subset\sR$, $\R\not=\sR$,  by the above corollary, 
since $\R$  is an infinite set. Similarly,  $\sN$, $^\ast\Bbb Z$, $^\ast\Bbb Q$, etc., 
are proper extensions of  $\N$, $\Bbb Z$, $\Bbb Q$, respectively. 

\vskip 0,5 truecm\noindent
{\bf (4.10) Theorem} {\em (Cartesian Products)}:

\vskip 0,2 truecm
{\bf (i)} The extension mapping  $^\ast$  preserves the Cartesian product, i.e.  
if  $A, B\in V(S) \setminus S$,  then  
$$
  ^\ast(A\times B) = \sA \times \sB.
$$
Consequently, the set of standard sets  $^\sigma V(S)\setminus S$  is closed under the 
Cartesian product of finite many sets.

\vskip 0,2 truecm
{\bf (ii)}  The extension mapping preserves the ordered pairing of  entities (individuals or 
sets), i.e. if  $a,b\in V(S)$,  then  
$$
  ^\ast\langle a,b\rangle = \langle\,^\ast a,\,^\ast b\rangle\;.
$$
Consequently, the set of standard sets  $^\sigma V(S)$  is closed under the building of 
ordered  $n$-tuples for  $n\in\N$. 

\vskip 0,3 truecm\noindent
{\bf Proof}: (i) Assume that  $A\times B = C$  which can be formalized in  $\cL(V(S))$  as 
$$
  [(\forall a\in A)(\forall b\in B)(\langle a,b\rangle\in C)]\wedge
  [(\forall c\in C)(\exists a\in A)(\exists b\in B)(\langle a,b\rangle = c)].
$$
Thus, 
$$
  [(\forall a\in\!\sA)(\forall b\in\!\sB)(\langle a,b\rangle\in\!\sC)]\wedge
  [(\forall c\in\!\sC)(\exists a\in\!\sA)(\exists b\in\!\sB)(\langle a,b\rangle = c)]
$$
holds in $\cL(V(\sS))$, by Transfer Principle, which means nothing but  
$\sA\times \sB = \sC$.  The generalization for  $n$  many sets follows by induction.

\vskip 0,2 truecm
(ii)$\;^\ast\langle a,b\rangle =\,^\ast\{\{a\},\{a,b\}\} = \{^\ast\{a\},^\ast\{a,b\}\} = 
    \{\{\sa\}, \{\sa,^\ast b\}\} = \langle\sa,^\ast b\rangle$, 
as required, by Theorem 4.6. $\tr$

\vskip 0,3 truecm\noindent
{\bf (4.11) Notation}:  Based on the above result, we have  $^\ast(A^n) = (\sA)^n$. 
So, we shall simply write  $\sA^n$ instead of  $^\ast(A^n)$ or $(\sA)^n$.
In particular for  $S=A=\R$,  and  $d\in\N$,  we write  $\sR^d$ instead of  $^\ast(\R^d)$
or  $(\sR)^d$.

\vskip 0,2 truecm\noindent
The next result is an addition to the Extension Principle. 

\vskip 0,3 truecm\noindent
{\bf (4.12) Lemma} {\em (Complex Numbers)}: Let  $S = \R$ and  $V(\sR)$  be a nonstandard 
model of  $\R$. Then  $^\ast z = z$  for all  $z\in\C$. 

\vskip 0,3 truecm\noindent
{\bf Proof}: We have $\C\in V(\R)$  since  $\C=\R^2$. Thus, both  $^\ast \C$  and  $^\ast z$ 
are well defined in  $V(\sR)$. Also, we have  $z=\langle x,y\rangle$  for some  $x,y\in\R$.
Thus, with the help of the above theorem, we have:
$$
  ^\ast z =\,^\ast\langle x,y\rangle = \langle^\ast x,^\ast y\rangle = \langle x,y\rangle = z,
$$
as required, since $^\ast x = x$ and $^\ast y = y$, by the {\em Extension Principle}.
$\tr$

\vskip 0,2 truecm   
Our next topic is some properties of the standard functions, i.e. the nonstandard 
extension of functions in  $V(S)$.

\vskip 0,3 truecm\noindent
{\bf (4.13) Theorem}: Let  $f : A\to B$  be a function in  $V(S)$,  i.e.  $A,B\in V(S)$.
Let  $\ssf$  be the nonstandard extension of  $f$. Then: 

\vskip 0,2 truecm
{\bf (i)} $\quad\ssf\;$  is a function of the type  $\ssf : \sA\to\sB$.

\vskip 0,2 truecm
{\bf(ii)} $\;\,\ssf\;$  is an extension of  $f$  in the sense that  $\ssf\mid\,^\sigma A = f,\quad$ i.e. 
$$
  \ssf(\sa) = \,^\ast (f(a)),
$$
for all  $a\in A$.

\vskip 0,2 truecm
{\bf (iii)}  Let  ${\it dom}(f)$ and  ${\it ran}(f)$  be the domain and the range of  $f$,
respectively, and let ${\it dom}(\ssf)$  and  ${\it ran}(\ssf)$  be the domain and the range 
of  $\ssf$, respectively. Then 
$$
  ^\ast({\it dom}(f)) = {\it dom}(\ssf)\qquad {\text{and}}\qquad^\ast({\it ran}(f))={\it ran}(\ssf).
$$

\vskip 0,3 truecm\noindent
{\bf Proof}:  (i) The fact that  $f$  is a function in  $V(S)$  and that  ${\it dom}(f)$
and  ${\it ran}(f)$  are its domain and range, respectively, can be formalized by the formula:
\begin{align*}
	&(\forall z \in f)(\exists x \in A)(\exists y \in B)[z = \langle x,y\rangle]\wedge \\
        &(\forall x \in A)(\exists y \in B)[\langle x,y \rangle \in f]\wedge \\
        &(\forall x \in A)(\forall y \in B)[(\langle x,y \rangle\in f)\Leftrightarrow (y=f(x))]
\end{align*}  
which is true in  $\cL(V(S))$. The first line of the above formula simply says that ``$f$  is a 
relation between $A$ and $B$", the second line says that ``$A$  is the domain of  $f$", the third 
line expresses the ``uniqueness of the value  $y = f(x)$  for any  $x$  in  $A$ ". 
By {\em Transfer Principle}, 
\begin{align*}
       &(\forall z\in\ssf)(\exists x\in\sA)(\exists y\in\sB) [z = \langle x,y\rangle]\wedge\\
       &(\forall x\in \sA)(\exists y\in\sB)[\langle x,y\rangle\in\ssf]\wedge\\
       &(\forall x\in \sA)(\forall y\in\sB)[(\langle x,y\rangle\in\ssf)\Leftrightarrow (y=\ssf(x))] 
\end{align*}
is true in $\cL(V(\sS))$. The above formula means nothing but that  $\ssf$  is a function of 
the type  $\ssf : \sA\to\sB$.

\vskip 0,3 truecm
(ii) Suppose  $a\in A$  and  $b\in B$. With the help of the Transfer Principle, we have
\begin{align*}
      &[ f(a)=b]\Leftrightarrow [(a\in A)\wedge (\langle a,b\rangle \in f)]\Leftrightarrow\\
      &\Leftrightarrow [(\sa\in\sA)\wedge (\langle\sa, ^\ast b\rangle\in\ssf)]\Leftrightarrow
        [\ssf(\sa) = ^\ast b],
\end{align*}
Hence,  $^\ast(f(a)) =\, ^\ast b = \ssf(\sa)$, as required.

\vskip 0,3 truecm
(iii) $^\ast({\it dom}(f)) = {\it dom}(\ssf)$  follows immediately from (i) since 
      ${\it dom}(f) = A$. Observe that  ${\it ran}(f)$  is described by
$$
   {\it ran}(f) = \{y\in B : (\exists x\in{\it dom}(f)) [\langle x,y\rangle\in f]\}.
$$
Hence, it follows
$$
   ^\ast({\it ran}(f)) = \{y\in\sB : (\exists x\in\,^\ast{\it dom}(f))[\langle x,y\rangle\in\ssf]\},
$$
by Lemma 4.4. Replacing  $^\ast({\it dom}(f)) = {\it dom}(\ssf)$, we get:
$$
  ^\ast({\it ran}(f)) = \{y\in\sB : (\exists x\in {\it dom}(\ssf))[\langle x,y\rangle\in\ssf]\}
$$
The latter formula means nothing but that  $^\ast({\it ran}(f)) = {\it ran}(\ssf$), as required.
$\qquad \blacktriangle$

\vskip 0,4 truecm
{\bf(4.16) Corollary} {\em (Functions in $S$)}: Let $f : A\to B$  be a function in the set 
of the individuals $S$,  i.e.  $A, B\subseteq S$. Then  $\ssf$  is an extension of  $f$  
in the usual sense, i.e.  $\ssf\mid A = f$, or
$$
  \ssf(a) =  f(a),
$$
for all  $a\in A$.

\vskip 0,3 truecm\noindent
{\bf Proof}: The result follows from Theorem 4.12 since  $\sa = a$  and  $^\ast(f(a))=f(a)$ 
for all  $a\in A$, by the {\em Extension Principle}. $\tr$

\vskip 0,5 truecm\noindent
{\bf § 5. Nonstandard Real Numbers}

\vskip 0,3 truecm
Let $S=\R$,  $V(\R)$ be its superstructure and  $\cL(V(\R))$  be its language. We shall refer 
to  $V(\R)$  as {\em Standard Analsyis}.  Let  $V(\sR)$  be a nonstandard extension of $V(\R)$
in the sense of Definition (4.1) and  $\cL(V(\sR))$  be its language. We shall refer to $V(\sR)$
as {\em Non-Standard Analsyis}.  Also the elements of  $\sR$  as {\em nonstandard real numbers}  
or {\em hyperreal numbers}. Similarly, $\sN, ^\ast\Bbb Z, ^\ast\Bbb Q$  denote the nonstandard 
extensions of  $\N,\Bbb Z, \Bbb Q$  respectively. We call their elements {\em nonstandard natural, 
nonstandard integer and nonstandard rational numbers}, respectively.  

\vskip 0,3 truecm
Let $A:\R\times\R\to\R$, $A(x,y) = x+y$, and  $M :\R\times\R\to\R, M(x, y) = x y$, be the 
addition and the multiplication in $\R$, respectively. Let $\R_+$  be the set of the positive 
real numbers. Let  $\sA,\sM$  and  $\sR_+$  be the nonstandard extensions of  $A,M$ and $\sR_+$,  
respectively. Observe that  $\sA$  and  $\sM$  are functions of the type  $\sA:\sR\times\sR\to\sR$ 
and  $\sM:\sR\times\sR\to\sR$, respectively, by Theorem 8.4 and Theorem 9.1.
 
\vskip 0,3 truecm\noindent
{\bf (5.1) Definition} {\em (Field Operations and Order Relation in $\sR$)}: We define the addition 
and multiplication in $\sR$, by  $x+y =\sA(x,y)$  and  $x\cdot y = M(x,y)$, respectively. 
The order relation in  $\sR$  is defined by  $x > 0$  if  $x\in\sR_+$. 

\vskip 0,3 truecm\noindent
{\bf (5.2) Theorem} {\em (Properties of $\sR$)}: The set of nonstandard real numbers  $\sR$ 
is a totally ordered non-Archimedean field which is a proper extension of $\R$, in symbols,  
$\R\subset\sR,\;\R\not=\sR$.

\vskip 0,3 truecm\noindent
{\bf Proof}: Let  $0$  and $1$ are the zero and the unit in $\R$, respectively. The fact that  
$\R$ is a totally ordered field can be formalized in  $\cL(V(\R))$  by the following statements:
\begin{itemize} 
  \item[] $(\forall x\in\R)([(x+0 = x)\wedge(x\; 0=0)]$
  \item[] $(\forall x\in\R)(\exists y\in\R)[A(x,y) = 0)]$
  \item[] $(\forall x\in\R)[M(x,1) = x]$
  \item[] $(\forall x\in\R)[(x\not= 0)\Rightarrow (\exists y\in\R)[M(x,y) = 1]]$
  \item[] $(\forall x\in\R)(\forall y\in\R)[A(x,y) = A(y,x)]$
  \item[] $(\forall x\in\R)(\forall y\in\R)[A(A(x,y),z) = A(x,A(y,z))]$
  \item[] $(\forall x\in\R)(\forall y\in\R)[M(x,y) = M(y,x)]$
  \item[] $(\forall x\in\R)(\forall y\in\R)[M(M(x,y),z) = M(x,M(y,z))]$
  \item[] $(\forall x\in\R)(\forall y\in\R)(\forall z\in\R)[M(A(x,y),z) = A(M(x,z), M(y,z))]$
  \item[] $0\in\R_+$
  \item[] $(\forall x\in\R_+)(\forall y\in\R_+)[A(x,y)\in\R_+)\wedge M(x,y)\in\R_+)]$
  \item[] $(\forall y\in\R)[(y=0)\vee (y \in \R_+) \vee (- y\in\R_+)],$
\end{itemize}
where  $- y$  is the (unique) solution of the equation  $A(x,y) = 0$ in  $\R$.
By {\em Transfer Principle}, it follows:
\begin{itemize}
  \item[] $(\forall x\in\sR)[(x+0=x)\wedge (x\,0=0)]$
  \item[] $(\forall x\in\sR)(\exists y\in\sR)[\sA(x,y) = 0]$
  \item[] $(\forall x\in\sR)[\sM(x,1) = x]$
  \item[] $(\forall x\in\sR)[(x\not= 0)\Rightarrow (\exists y\in\sR)[\sM(x,y) = 1]]$
  \item[] $(\forall x\in\sR)(\forall y\in\sR)[\sA(x,y) = \sA(y,x)]$
  \item[] $(\forall x\in\sR)(\forall y\in\sR)[\sA(\sA(x,y),z) =\sA(x,\sA(y,z))]$
  \item[] $(\forall x\in\sR)(\forall y\in\sR)[M(x,y) = M(y,z)]$
  \item[] $(\forall x\in\sR)(\forall y\in\sR)[\sM(\sM(x,y),z) = \sM(x,\sM(y,z))]$
  \item[] $(\forall x\in\sR)(\forall y\in\sR)(\forall z\in\sR)[\sM(\sA(x,y),z) = $
      \[= \sA(\sM(x,z),\sM(y,z))]\]
  \item[] $0\notin\sR_+$
  \item[] $(\forall x\in\sR_+)(\forall y\in\sR_+)[(\sA(x,y)\in\sR_+)\wedge (\sM(x,y)\in\sR_+)]$
  \item[] $(\forall y\in\sR)[(y\not= 0)\vee (y\in\sR_+)\vee (- y\in\sR_+),$
\end{itemize}
		
where  $- y$  is the (unique) solution of the equation  $\sA(x, y) = 0$  in  $\sR$.
The interpretation of the above formulae mean nothing but that  $\sR$  is a totally ordered field.  
On the other hand, $\R\subset\sR$, $\R\not=\sR$  follows from Corollary 4.8 (applied for  
$A = S = \R$), since  $\R$  is an infinite set. Thus,  $\sR$  turns out to be a proper 
totally ordered field extension of  $\R$. It follows that  $\sR$  is a non-Archimedean field 
(any proper totally ordered field extension of  $\R$  is  non-Archimedean).$\tr$

\vskip 0,3 truecm
Let  $\cI(\sR),\; \cF(\sR)$  and  $\cL(\sR)$  denote, as usual, the sets of the infinitesimals, 
finite and infinitely large numbers in  $\sR$, respectively. Recall that  $\alpha\in\cI(\sR)$ 
if  $\mid\alpha\mid<1/n$  for all  $n\in\N$,  $\alpha\in\cF(\sR)$  if  $\mid\alpha\mid< n$  
for some  $n\in\N$, and  $\alpha\in\cL(\sR)$  if  $\mid\alpha\mid > n$  for all $n\in\N$. 
The {\em infinitesimal relation}  in  $\sR$  is defined by:  If  $\alpha,\,\beta\in\sR$, then 
$\alpha\approx\beta$ if  $\alpha - \beta\in\cI(\sR)$. Notice that (as in any totally 
ordered field) we have 
\begin{align*}
  &\sR = \cF(\sR)\cup\cL(\sR),\quad\cF(\sR)\cap\cL(\sR) = \emptyset,\\
  &\cI(\sR)\subset\cF(\sR),\quad \R\subset\cF(\sR),\\
  &\R\cap\cI(\sR) = \{0\},\\
  &\cL(\sR) = \{1/x : x\in\cI(\sR),\; x\not= 0\}.
\end{align*}
The fact that  $\sR$  is a non-Archimedean field means that  $\sR$  has non-zero infinitesimals 
and infinitely large elements, in symbols, $\cI(\sR)\setminus\{0\}\not=\emptyset$  and 
$\cL(\sR)\not=\emptyset$. Recall also that (as in any totally ordered field),  $\cF(\sR)$ is a 
convex Archimedean integral domain (totally ordered Archimedean ring without zero divisors) 
and  $\cI(\sR)$  is a convex maximal ideal in  $\cF(\sR)$. Hence, the factor space 
$\cF(\sR)/\cI(\sR)$  is a totally ordered Archimedean field. Recall further that  (as in any 
totally ordered field) we have 
$$
  \{r + h : r\in\R,\; h\in\cI(\sR)\} \subseteq \cF(\sR),
$$
and
$$
  \quad\left\{\frac{a+h}{b+g}\,:\,a,\, b\in\R,\, h,\, g\in\cI(\sR)\right\}\subseteq\sR.
$$
Observe that  $\alpha\in\cF(\sR)$ in  $\alpha = r + h$  determines uniquely  $r\in\R$ and 
$h\in\cI(\sR)$, due to (5.5). In addition, the order completeness of  $\R$  implies that the 
inclusions in (5.7) and (5.8) are, actually,  equalities:

\vskip 0,3 truecm\noindent
{\bf (5.9) Theorem}: We have the following characterizations of  $\cF(\sR)$  and  $\sR$:

\vskip 0,2 truecm
   {\bf (i)} $\quad \cF(\sR) = \{a+h : a\in\R,\, h\in\cI(\sR)\}$.  

\vskip 0,2 truecm
  {\bf (ii)} $\; \sR = \left\{\dfrac{a+h}{b+g}\,:\; a,\, b\in\R,\, h,\, g\in\cI(\sR)\right\}$.

\vskip 0,3 truecm\noindent
{\bf Proof}: (i) Suppose  $\alpha\in\cF(\sR)$. We have to show that  $ \alpha = a+h$ 
for some $a\in\R$, $h\in\cI(\sR)$. Let  $a = {\it sup}\{x\in\R : x <\alpha\}$ and  
$h = \alpha - a$. Notice the order completeness of  $\R$ guaranties the existence of  $a$.
It suffices  to show that  $h\in\cI(\sR)$. Suppose (for contradiction) that  
$h\notin\cI(\sR)$, i.e.
there exists  $\eps\in\R_+$  such that  $\eps < |\alpha - a |$. If  $\alpha - a > 0$, 
then we have $a + \eps < \alpha$ contradicting the fact that  $a$  is an upper bound of the set 
$\{x\in\R : x < \alpha\}$.  If  $\alpha - a < 0$, then we have   $\alpha < a - \eps$,  
contradicting the maximality of  $a$. 

\vskip 0,2 truecm
(ii)  follows immediately from (i) and $\sR = \cF(\sR)\cup\cL(\sR)$. Indeed, suppose that 
$\alpha\in\sR$. If $\alpha$  is a finite number, then  $\alpha = a + h$, by (i), thus, 
$\alpha = \dfrac{a+h}{b+g}\;$  for  $b=1$ and $g=0$. If $\alpha$ is an infinitely large number, then  
\vskip 0,2 truecm\noindent
$\alpha = \dfrac{a+h}{b+g}\;$  for  $a=1,\, h=0,\, b=0$  and  $g=1/\alpha$. 

\vskip 0,5 truecm\noindent
{\bf (5.10) Definition} {\em (Standard Part)}: We define the standard part mapping  
$\text{st} : \sR\to\R\cup\{\pm\infty\}$  by  $\text{st}(a+h)=a$  if  $a\in\R$  and  $h\in\cI(\sR)$ 
and by  $\text{st}(\alpha)=\pm\infty$  if  $\alpha\in\cL(\sR),\,\alpha>0$  or  $\alpha<0$, 
respectively. 

\vskip 0,5 truecm\noindent
{\bf (5.11) Lemma:}

\vskip 0,2 truecm
{\bf (i)} $\;\alpha\in\cF(\sR)$  iff  $\text{st}(\alpha)\in\R$  and in this case we have the (unique) 
presentation:
$$
  \alpha = \text{st}(\alpha) + h
$$
for some  $h\in\cI(\sR)$. Or, equivalently, every finite number  $\alpha\in\cF(\sR)$ is 
infinitely close to a unique real number  $\text{st}(\alpha)$, in symbols,  $\alpha\approx \text{st}(\alpha)$. 

\vskip 0,3 truecm
{\bf (ii)} The totally ordered field  $\cF(\sR)/\cI(\sR)$  is isomorphic to $\R$  under the 
mapping  $q(\alpha)\to \text{st}(\alpha)$, where  $q : \cF(\sR)\to\cF(\sR)/\cI(\sR)$  is the 
corresponding quotient mapping.

\vskip 0,5 truecm\noindent
{\bf Proof}: Both (i) and (ii) are simple reformulatings of the previous result taking into 
account that  $\R\subset\sR$.   $\tr$

\vskip 0,3 truecm\noindent
{\bf (5.13) Theorem} {\em (Properties of $\text{st}$)}: Let  $\alpha,\,\beta\in\cF(\sR)$.  Then we have: 

\vskip 0,2truecm
{\bf (i)} $\quad \alpha\approx\beta$ iff  $\text{st}(\alpha)=\text{st}(\beta)$. In particular,  
          $\alpha\in\cI(\sR)$  iff  $\alpha\approx 0$  in  $\sR$  iff   
          $\text{st}(\alpha)=0$  in  $\R$.

\vskip 0,2 truecm
{\bf (ii)} $\;\, \text{st}(\alpha\pm\beta) = \text{st}(\alpha)\pm \text{st}(\beta)$.

\vskip 0,2 truecm
{\bf (iii)} $\; \text{st}(\alpha\beta) = \text{st}(\alpha) \text{st}(\beta)$;

\vskip 0,2 truecm 
{\bf (iv)} $\; \text{st}(\alpha/\beta) = \text{st}(\alpha)/\text{st}(\beta)$  whenever  $\text{st}(\beta)\not= 0$. 

\vskip 0,2 truecm
{\bf (v)} $\;\, \text{st}(\alpha^n) = (\text{st}(\alpha))^n$  for all  $n\in\N$.

\vskip 0,2 truecm 
{\bf (vi)} $\; \text{st}(\root n\of{\alpha}) = \root n\of{\text{st}(\alpha)},\; n\in\N$, whenever  
           $\root n\of\alpha$  exists in  $\sR$. In more details, if  $n$  is odd, 
           then the above equality holds for all  $\alpha\in\sR$, while the condition  
           $\text{st}(\alpha)>0$  is required in the case of even  $n$. 

\vskip 0,2 truecm
{\bf (vii)} $\,$ If $\alpha\not\approx\beta$, then  $\alpha < \beta$  iff  $\text{st}(\alpha)< \text{st}(\beta)$.
            As a result,  $\alpha\le\beta$  in  $\sR$  implies  $\text{st}(\alpha)\le \text{st}(\beta)$
            in  $\R$.

\vskip 0,3 truecm\noindent 
{\bf Proof}: The properties (i)-(iii) follows immediately from the definition of $\text{st}$. 
To show (iv), apply  st  to both sides of  $\alpha = \beta(\alpha/\beta)$. It follows 
$\text{st}(\alpha) = \text{st}(\beta)\text{st}(\alpha/\beta)$, by (iii), which implies (iv). The property (v) 
follows from (iii) by induction. To show (vi), notice that  $\beta = \root n\of\alpha$  is 
equivalent to  $\beta^n = \alpha$. Thus, applying (v), we have  $(\text{st}(\beta))^n = \text{st}(\alpha)$,
which is equivalent to  $\text{st}(\beta)=\root n\of{\text{st}(\alpha)}$, as required. Finally, (vii) 
follows directly from the convexity of  $\cI(\sR)$. $\tr$

\vskip 0,8 truecm\noindent
{\bf (5.14) Example} {\em (Functions)}: Let  $f:\R_+\to\R$  be defined by  $f(x)=\ln(x)$.
For the nonstandard extension we have $\ssf:\sR_+\to\sR$ is defined by  $\ssf(x)=\,^\ast\ln(x)$ 
(Theorem 4.12). In other words,  $^\ast\ln(x)$  is well defined on  $\sR_+$  and for any 
$y\in\sR$  the equation  $y =\, ^\ast\ln(x)$  has a (unique) solution  $x$  in  $\sR_+$. 
In particular $^\ast\ln(x)$  is well defined for all positive infinitesimals in  $\sR$ 
(and the value of  $^\ast\ln(x)$  is a negative infinitely large number). Finally, 
$^\ast\ln$ is an extension of  $\ln$, i.e.  $^\ast\ln(x) = \ln(x)$  for all  
$x\in\R_+$  (Corollary 4.15).

\vskip 1 truecm
\noindent
{\bf{PART II. NONSTANDARD AND STANDARD COMPACTIFICATIONS OF TOPOLOGICAL SPACES}}

\vskip 0,5 truecm            
We use the nonstandard methods to construct all Hausdorff compactifications of a 
given topological space  $(X,T)$. This result is a generalization of an earlier work by K.D. 
Stroyan [21] about the compactifications of completely regular spaces. We also describe the 
Hewitt realcompactification which seems to be treated here for the first time in the 
nonstandard literature.

\vskip 0,2 truecm
There are a vast nonstandard works done on the Hausdorff compactifications: A. 
Robinson [17]-[18], W.A.J. Luxemburg [14], M. Machover and J. Hirschfeld [15], K.D. 
Stroyan [21], K.D. Stroyan and W.A.J. Luxemburg [22], H. Gonshor [4] and L. Haddad [5] 
and others. We believe that our description of the Hausdorff compactifications, in particular, 
the Stone-\v Cech compactifications of  $(X,T)$  is noticeably  simpler than those both in the 
standard and nonstandard literature mostly due to the fact that we manage to avoid involving 
the weak topology both on the initial space and its compactification.

\vskip 0,2 truecm 
Our technique can be shortly described as follows: To any topological space $(X,T)$  we 
attach its {\em nonstandard compactification} $(\sX, ^sT)$, where  $\sX$  is the nonstandard 
extension  of  $X$  supplied with the {\em standard topology}  $^sT$,  generated by all sets 
of the form  $\sG$,  where  $G\in T$. The standard topology  $^sT$  is courser than the 
{\em discrete} $S$-{\em topology} on $\sX$, (known also as LS-{\em topology}, where  $L$ 
stands for {\em Luxemburg}) with basic open sets:  $^\sigma\cP(X) = \{\sS:S\in\cP(X)\}$, 
introduced by W.A.J. Luxemburg ([14], p.47 and p.55) for a similar purpose. Our space 
$(\sX, ^sT)$  is compact (non Hausdorff) and every continuous function  $f$  on  $(X,T)$ 
has a unique continuous extension  $\sf$  on  $(\sX, ^sT)$. In contrast to the case of the 
discrete  $S$-topology, however,  $(\sX, ^sT)$  contains  $(X,T)$  densely. These properties 
of  $(\sX, ^sT)$ simplify essentially our next steps: We supply the nonstandard hull  
$\widehat X_\Phi = \widetilde X_\Phi / \sim_\Phi$   with the quotient topology   $\widehat T$,  
and show that  the space  $(\widehat X_\Phi, \widehat T)$  is Hausdorff.  Here,  the set  
of  the $\Phi$ - {\em finite points}  $\widetilde X_\Phi \subseteq\sX$  and the equivalence 
relation ``$\sim_\Phi$" are specified by a  family of continuous functions  $\Phi$  and, 
thus, changing  $\Phi\subseteq C_b(X,\R)$, we describe in a uniform way all Hausdorff 
compactifications of $(X,T)$ as well as the Hewitt realcompactification of  $(X,T)$. 
If  $\Phi$  consists of bounded functions only, we have   $\widetilde X_\Phi=\sX$ and 
$\widehat X_\Phi=q[\sX]$, thus, the compactness of  $\widehat X_\Phi$  follows simply 
with the argument that the {\em continuous image of a compact space is compact}. 
In particular, when  $\Phi = C_b(X,\R)$  we obtain the Stone-\v Cech compactification 
$\beta(X,T)$  of  $(X,T)$  and when  $\Phi = C(X,\R)$  we obtain the Hewitt real 
compactification   $\nu(X,T)$  of  $(X,T)$. 

\vskip 0,2 truecm
We should mention that a technique based on the nonstandard compactification  $(\sX, ^sT)$
of  $(X,T)$  has already been successfully exploited for studying the compactifications of 
{\em ordered topological spaces}  by the authors of this paper (S. Salbany, T. Todorov [19]-[20]). 

\vskip 0,2 truecm
We shall use as well the terminology of  (J.L. Kelley [12]) and (L. Gillman and M. 
Jerison [3]). For the connection between the standard and nonstandard methods in topology 
we refer to (L. Haddad [5]).

\vskip 1 truecm\noindent

{\bf 1. Preliminaries: Monads and Their Basic Properties.}

\vskip 0,3 truecm
We shall briefly recall the definition of monads and some of their properties. For the 
original sources we refer to (A. Robinson [17]) and (K.D. Stroyan and W.A.J. Luxemburg 
[22], Chapter 8). For the general theory of monads, we refer to  W.A.J. Luxemburg [14] and 
K.D. Stroyan [21].

\vskip 0,2 truecm
Let  $(X,T)$  be a topological space. In order to apply nonstandard methods we need the 
superstructure  $V(S)$  over some set $S$  such that  $S = X\cup\R$  (the choice 
$S\supset X\cup\R$  also will do), and a $\kappa$-saturated nonstandard model  
$V(\sS)$  of  $S$  with  $\kappa > \text{card}\,T$ (Chapter I, Section 2). Sometimes we shall 
consider two topological  spaces  $(X,T)$  and  $(X',T')$. In this case we shall assume that  
$S = X\cup X' \cup\R$  (or  $S\supseteq X\cup X'\cup \R$)  and  
$$
  \kappa > \text {max}(\text{card}\,T,\, \text{card}\,T'). 
$$
Any polysaturated model will cover all those cases (Chapter I, Definition 2.8). We shall often 
refer to the {\em Extension, Transfer and Saturation Principles} (Chapter I, Section 2 as Axiom 
1-3, respectively), and also the {\em Boolean Properties} of the extension mapping (Theorem I.4.2).

\vskip 0,3 truecm\noindent 
{\bf (1.1) Definition} {\em (Monads)} : Let  $(X,T)$  be a topological space and  $\sX$  be the 
nonstandard extension of  $X$. Then:

\vskip 0,2 truecm 
{\bf (i)}  For any  $\alpha\in\sX$   define  the monad  $\mu(\alpha)$  of  $\alpha$  by
$$
  \mu(\alpha) = \bigcap\{\sG\mid \alpha\in\sG,\,  G\in T\}. \leqno{\bf (1.2)}
$$

\vskip 0,2 truecm
{\bf (ii)}  For any  $A\subseteq\sX$   define
$$
  \mu(A) = \bigcap \{\sG\mid A\subseteq G,\,  G\in T\}. \leqno{\bf (1.3)}
$$

The monad of a set  $A$  is obviously a generalization of the monad at a point  $\alpha$  when  
$A = \{\alpha\}$  for some  $\alpha\in\sX$. We use the same notation for both. Also for any  
$A\subseteq X$  we have
$$
  \mu(A) = \mu(\sA). \leqno{\bf (1.4)}
$$
The following properties of  monads follow almost directly from the definition. 

\vskip 0,3 truecm\noindent 
{\bf (1.5) Lemma}: If  $A,B\subseteq\sX$, then:

\vskip 0,2 truecm
{\bf (i)}  $\;\; A\subseteq\mu(A)$.

\vskip 0,2 truecm
{\bf (ii)}  $\;A\subseteq B$  implies  $\mu(A)\subseteq\mu(B)$.

\vskip 0,2 truecm
{\bf (iii)}  $\,\mu(\mu(A)) = \mu(A)$.

\vskip 0,3 truecm
The above lemma shows that the monad of a set is a generalized closure operator in  $\sX$  
(see e.g. P.C. Hammer [6] and K.D. Stroyan [21], Section 2).

\vskip 0,3 truecm\noindent 
{\bf (1.6) Corollary}:  For any  $A\subseteq\sX$  and any  $\alpha, \beta\in\sX$:

\vskip 0,2 truecm
{\bf(i)} $\;\;\alpha\in A$  implies   $\mu(\alpha)\subseteq \mu(A)$.

\vskip 0,2 truecm
{\bf (ii)} $\; \alpha\in\mu(\beta)$  iff  $\mu(\alpha)\subseteq \mu(\beta)$.

\vskip 0,2 truecm
{\bf (iii)}  $\, \alpha\in\mu(\beta)$  and  $\beta\in\mu(\alpha)$  iff  $\mu(\alpha) = \mu(\beta)$.

\vskip 0,3 truecm\noindent
{\bf Proof}:  (i)  follows from (1.5)-(ii) by  $A = \{\alpha\}$;  (ii)  follows from (1.5)-(ii) 
and (1.5)-(iii). Indeed,  $\alpha\in\mu(\beta)$  implies  $\{\alpha\}\subseteq\mu(\beta)$  
which implies  $\mu(\alpha)\subseteq\mu(\mu(\beta))= \mu(\beta)$. The converse is clear; 
(iii) follows directly from (ii).  $\tr$

\vskip 0,3 truecm\noindent 
{\bf (1.7) Theorem} {\em (Balloon and Nuclei Principles)}: Let $(X,T)$  be a topological space, 
$x\in X$, and  $\mu(x)$  be the monad of  $x$  at  $(X,T)$. 

\vskip 0,2 truecm
{\bf (i) } {\em Balloon Principle} :  If  $\mu(x)\subset\cB$  for some internal set  
$\cB\subseteq\sX$, then there exists  $G\in T$  such that  $\mu(x)\subset\sG\subseteq\cB$
(ballooning of   $\mu(x)$  into  $\sG$).

\vskip 0,2 truecm
{\bf (ii)} {\em Nuclei Principle} : There exists an internal set  $\cA\subseteq\sX$  such that  
$x\in\cA\subset\mu(x)$. The set  $\cA$  is called a nuclei  of  $\mu(x)$. 

\vskip 0,3 truecm\noindent
{\bf Proof}: (i) Suppose not, i.e.  $\sG - \cB\not= \emptyset$  for all  $G\in T, x\in G$. 
Observe that the family of sets  $\{\sG - \cB\}_{G\in T,\,x\in G}\,,$ has the finite intersection 
property since  $(\sG_1 - \cB)\cap (\sG_2 - \cB) =\, ^\ast(G_1\cap G_2)-\cB$.  It follows
$$
  \mu(x)-\cB = \bigcap\limits_{x\in G\in T} (\sG - \cB)\not=\emptyset,
$$
by {Saturation Principle}, since 
$$
  \text{card}\{G : x\in G\in T\}\le \text{card}\,T\le\kappa, 
$$
by the choice of the nonstandard model. But  $\mu(x) - \cB\not=\emptyset$  contradicts our 
assumption.         

\vskip 0,2 truecm
(ii) Define the family  $\{S_G\}_{x\in G\in T}$, where  $S_G = \{H\in T : x\in H\in G\}$,
and observe that it has the finite intersection property since  $G\in S_G$, thus, 
$S_G\not=\emptyset$, and, on the other hand,  $S_{G_1}\cap\, S_{G_2} = S_{G_1\cap G_2}$. 
It follows that there exists  $\cA$  in the intersection
$$
   \bigcap\limits_{x\in G\in T} \sS_G,
$$
by the {\em Saturation Principle}.  On the other hand, observe that 
$$ 
  \sS_G = \{H\in\sT : x\in H\subseteq\sG\}.
$$
Thus, $\cA$  is internal (as an element of  $\sT$) and  $\cA\subset\mu(x)$, as required. 
$\tr$

\vskip 0,5 truecm\noindent
The next result is due to  A. Robinson ([17], Theorem 4.14., p.90):

\vskip 0,3 truecm\noindent
{\bf (1.8) Theorem} {\em (A. Robinson)}: 

\vskip 0,3 truecm
{\bf (i)}  Let  $(X,T)$ be a topological space and let  $x\in H\subseteq X$  and  $x\in X$. 
Then  $x$  is an interior point of  $H$  in $(X,T)$  iff  $\mu(x)\subset\sH$. Consequently, 
$H$  is open in  $(X,T)$  iff  $\mu(x)\subset\sH$  for all  $x\in H$.

\vskip 0,2 truecm
{\bf (ii)}  A set  $F\subseteq X$  is closed in  $(X,T)$  iff  $\sF\cap\mu(x)\not=\emptyset$  
implies  $x\in F$  for any  $x\in X$.

\vskip 0,2 truecm
{\bf (iii)}  Let  $A\subseteq X$  and  $\text{cl}_X(A)$  be the closure of $A$ in $(X,T)$. Then 
$$
   \text{cl}_X(A) = \{x\in X : \sA\cap \mu(x)\not=\emptyset\}. \leqno{\bf (1.9)}
$$

\vskip 0,3 truecm\noindent
{\bf Proof}: (i) $(\Rightarrow)$ If  $x$  is an interior point of  $H$, then $\mu(x)\subset\sH$, 
by the definition of  $\mu(x)$. 

\vskip 0,2 truecm
$(\Leftarrow)$ Suppose (for contradiction) that  $x$  is not an interior point of  H, i.e.  
$G - H\not=\emptyset$ for all  $G$  such that  $x\in G\in T$. Observe that the family of sets 
$\{G-H\}_{x\in G\in T}$  has the finite intersection property. It follows that the family of 
internal (actually, standard) sets $\{\sG-\sH\}_{x\in G\in T}$  has the finite intersection 
property, since  $^\ast(G-H) = \sG - \sH$, by the {\em Boolean Properties}. (Theorem I.4.2) 
It follows that its intersection  $\mu(x)-\sH$  is non-empty, by the {\em Saturtion Principle}, 
a contradiction. 

\vskip 0,2 truecm
(ii) Suppose (for contradiction) that  $x\in X-F$.  We have  $\mu(x)\subset\sX-\sF$, by the 
above theorem, since  $X-F$  is open, by assumption, and  $\sX-\sF =\, ^\ast(X-F)$, 
by Theorem I.4.2. It follows  $\mu(x)\cap\sF = \emptyset$, a contradiction. 

\vskip 0,2 truecm
(iii) $(\subseteq$) Let  $x\in\text{cl}_X(A)$,  i.e.  $x\in F$  for all  $F$  such that  
$A\subset F\subseteq X$, $X-F\in T$. Suppose (for contradiction) that  $\sA\cap\mu(x)=\emptyset$. 
Then, by the {\em Balloon Principle} (applied for  $\cB = \sX-\sA)$, there exists  $G\in T$,
$x\in G$, such that  $\sA\,\cap \sG = \emptyset$. Thus, we have  $\sA\subseteq\,^\ast(X-G)$,  
implying  $A\subseteq X-G$, by the {\em Boolean Properties}. Hence, it follows 
$x\in X-G$, by our assumption (since  $X-G$  is a closed set), a contradiction. 

\vskip 0,2 truecm
$(\supseteq)$  Let  $x\in X$  and  $\sA\cap \mu(x)\not=\emptyset$. We have to show that  
$x\in F$  for all  $F$  such that  $A\subset F\subseteq X$  and  $X-F\in T$. Suppose 
(for contradiction) that  $x\notin F$  for some  $F$  such that  $A\subset F\subseteq X$ 
and  $X-F\in T$.  It follows that  $x\in X-F$. On the other hand,  $A\subset F$ implies 
$\sA\subset\sF$, by the {\em Boolean Properties}.  Hence,  $\sA\cap(\sX-\sF) = \emptyset$, 
which implies  $\sA\cap\mu(x)\not=\emptyset$ (because  $\sA\cap\mu(x)\subseteq\sA\cap(\sX-\sF)$), 
a contradiction.  $\tr$

\vskip 0,5 truecm\noindent
{\bf (1.10) Definition} {\em (Nearstandard Points and Standard Part)}: Let  $(X,T)$  be a 
topological space and  $\mu(x)$, $x\in X$,  be its monads. 

\vskip 0,2 truecm
{\bf (i)} If  $A\subseteq X$, then the points in the union  $\widetilde A = \cup_{x\in A}\;\mu(x)$
are called {\em nearstandard points of}  $\sA$. In particular, the points in  
$\widetilde X = \cup_{x\in X}\;\mu(x)$  are called nearstandard points of  $\sX$. 

\vskip 0,2 truecm
{\bf (ii)} Assume, in addition, that  $(X,T)$  is a regular Hausdorff space. Then the mapping  
$\text{st}_X :\widetilde X\to X$, defined by   $\text{st}_X(\xi) = x,\,\xi\in\mu(x)$, 
is called {\em standard part mapping}. 

\vskip 0,2 truecm
Notice that the assumption that  $(X,T)$  is a {\em regular Hausdorff space} guarantees the 
correctness of  $\text{st}_X$  (the uniqueness of  $x$). We shall often skip the subindex  
and write simply  st  if no confusion could arise.  

\vskip 0,3 truecm\noindent
{\bf (1.11) Examples}: 

\vskip 0,3 truecm
{\bf 1.} Let $(\R,\tau)$ be the space of the real numbers supplied with the usual topology $\tau$.
Then the nearstandard points are, actually, the finite points, in symbols, $\widetilde\R=\cF(\sR)$
and  $\text{st}(\xi)=x$, for  $\xi\in\cF(\sR)$, $x\in\R$, $\xi \approx x$. 

\vskip 0,2 truecm
{\bf 2.} Let  $(I,\tau)$, where  $I=(a,b)=\{x\in\R:a<x<b\}$. Then the nearstandard points are
$$
  \widetilde I = \{x\in\sR : a<x<b,\; x\not\approx a,\; x\not\approx b\}  
$$
and, as before,  $\text{st}(\xi) = x$,  for  $\xi\in\widetilde I,\; x\in\R,\;\xi\approx x$. 

\vskip 0,2 truecm
{\bf 3.} This example illustrates the {\em Nuclei Principle}:  Let  $x\in I$  and  
$\mu(x) = \{\xi\in\sR : \xi\approx x\}$  be the monad of  $x$  in  $(I,\tau)$. Let  
$\rho\in\sR$,$\,\rho>0$,$\,\rho\approx 0$, be a positive infinitesimal, and observe that the set 
$$
  \cA = \{\xi\in\sR : \mid\xi-x\mid<\rho\}.
$$
is internal, by Theorem I.4.10. It follows that  $\cA$  is a nuclei of  $\mu(x)$ since (obviously) 
$x\in\cA\subset\mu(x)$.

\vskip 0,3 truecm\noindent
{\bf (1.12) Corollary}:  Let  $A\subset\R$  and let  $\,\text{cl}_{\Bbb R} (A)$  be the closure of  
$A$  in $(\R,\tau)$, where  $\tau$  is the usual topology of  $\R$. Then 

\vskip 0,2 truecm
{\bf (i)} $\;\text{cl}_{\Bbb R} (A)=\{x\in\R:\text{st}(\alpha)=x$  for some $\alpha\in\sA\}=$

\hspace{2cm}$ =\{x\in\R:x\approx\alpha$  for some  $\alpha\in\sA\}$.

\vskip 0,2 truecm
{\bf (ii)} If $A$  is bounded in  $\R$, then  $\text{cl}_{\Bbb R} (A)=\text{st}[\sA]$.

\vskip 0,5 truecm\noindent
{\bf Proof}: (i) $\text{cl}_{\Bbb R}(A) = \{x\in\R:\sA\,\cap\,\mu(x)\not=\emptyset\} =$

\hspace{2,5cm} $ = \{x\in\R:\text{st}(\alpha) = x$  for some  $\alpha\in\sA\}.$

\vskip 0,2 truecm 
(ii) There exists   $b\in\R$  for wich the formula
$$
  \Phi(A,b) = (\forall x\in A)(\,|x|\le b)
$$
is true in $\cL(V(\R))$. It follows that the formula
$$
  \Phi(\sA,b) = (\forall x\in\sA)(\,|x|\le b)
$$
is true in  $\cL(V(\sR))$, by {\em Transfer Principle}, since  $^\ast b = b$, by 
{\em Extension Principle}. The latter implies that $\sA\subset\cF(\sR)$, thus, $\text{st}[\sA]$
is well defined. On the other hand, 
\begin{align*}
     &\text{st}[\sA] = \{\text{st}(\alpha) :\quad \text{for some}\quad \alpha\in\sA\} = \\
     &= \{x\in\R : x = \text{st}(\alpha)\quad\text{for some}\quad \alpha\in\sA\} =
        \text{cl}_\Bbb R(A),
\end{align*}
by (i).  $\tr$

\vskip 0,8 truecm\noindent
{\bf (1.13) Lemma:}

\vskip 0,2 truecm  
{\bf (i)}   Let  $A,B\subseteq\sX$.  Then  $\mu(A)\cap\mu(B) = \emptyset$  iff  there exist 
open disjoint sets  $G$ and  $H$  such that  $A\subseteq\sG$  and  $B\subseteq\sH$.

\vskip 0,3 truecm
{\bf (ii)}  Let  $\alpha,\beta\in\sX$. Then  $\mu(\alpha)\cap\mu(\beta) = \emptyset$  
iff  there exist open disjoint sets  $G$ and  $H$  such that  $\alpha\in\sG$  and  
$\beta\in\sH$ .

\vskip 0,5 truecm\noindent
{\bf Proof}:  (i) Let  $\mu(A)\cap\mu(B) = \emptyset$  and suppose that  $G\cap H\not=\emptyset$  
for all open  $G$ and  $H$  such that  $A\subseteq\sG$  and  $B\subseteq\sH$. By the Saturation  
Principle (Chapter I, Section 2, Axiom 3),  we have
$$
   \mu(A)\cap\mu(B) = \cap\{^\ast(G\cap H) : G,H\in T,\, A\subseteq\sG\quad
   \text{and}\quad  B\subseteq\sH\}\not=\emptyset
$$
which is a contradiction.  The converse follows immediately;  

\vskip 0,2 truecm
(ii)  follows directly from (i) by    $\text{letting}\,  A = \{\alpha\}$  and  $B = \{\beta\}$.
$\tr$

\vskip 0,3 truecm
We shall have occasion to use other monads: Following (K.D. Stroyan and W.A.J.
Luxemburg [12], p. 195) , we state:

\vskip 0,3 truecm
{\bf (1.14)  Definition} {\em (General Monads)}:  Let  $X$  be a  set  and  $\cE$  be a ring 
of subsets of  $X$. Then  for  any   $\alpha\in\sX$   and  any   $A\subseteq\sX$  we  define 
the  $\cE$-monads of  $\alpha$  and  $A$,  respectively, by:
\begin{align*}
\mu_{\cal E}(\alpha)&=\,\bigcap\{\sG : G\in \cE\; , \;\alpha\in\sG\};,\\
\mu_{\cal E}(A)     &=\,\bigcap\{\sG : G\subseteq \cE\; , \;A\subseteq\sG\};.
\end{align*}
		
In the particular case of  $\cE = T$, where  $T$  is a  topology of  $X$, we obtain  $\mu_T=\mu$.   

As in the previous lemma we have:

\vskip 0,3 truecm\noindent
{\bf (1.15) Lemma}: Let  $X$   be a  set  and  $\cE$  be a ring of subsets of  $X$  and  
$A,B\subseteq\sX$.  Then  $\mu_{\cal E} (A)\cap\mu_{\cal E}(B) = \emptyset$  iff  there exist  
disjoint  $G,H\in\cE$  such that  $A\subseteq\sG$  and  $B\subseteq\sH$.\\

\noindent
{\bf 2. Nonstandard Compactification}

\vskip 0,3 truecm
By  $\sX$  and $\sR$  we denote the nonstandard extensions of  $X$  and  $\R$, respectively. 
If  $G\subseteq X$  and  $A\subseteq\R$, then  $\sG\subseteq\sX$  and  $\sA\subseteq\sR$  
will be the nonstandard extensions of  $G$  and  $A$, respectively (Definition I.2.4). 
For the more general concept of {\em internal set}  we refer again to (Definition I.2.4). 
Let  $(X,T)$  and  $(X',T')$  be two topological  spaces and $f : X\to X'$ be a function. 
Then  $\ssf : \sX\to\sX'$  will be the nonstandard extension of  $f$  (Theorem I.4.12).

\vskip 0,5 truecm\noindent
{\bf (2.1) Notations}: Let  $(X,T)$  be a topological space. Then, a simple observation shows 
that the collection of sets:  
$$
  ^\sigma T = \{\sG : G\in T\} \leqno{\bf (2.2)}
$$
forms a base for a topology in  $\sX$. We shall denote this topology  by  $\hsT$  and the 
corresponding topological  space by  $(\sX,\hsT)$.  Notice that the collection of sets:
$$
  \cF = \{\sF : X-F\in T\} \leqno{\bf (2.3)}
$$
forms a base of the closed sets of  $\sX$  in  $(\sX,\hsT)$.

\vskip 0,5 truecm\noindent
{\bf (2.4) Definition}{\em (Nonstandard Compactification)}:  Let  $(X,T)$   be a topological 
space and  $(\sX,\hsT)$  be the corresponding topological  space defined as above. Then:

\vskip 0,2 truecm
{\bf (i)}  $\;\hsT$  will be called the {\em standard topology on}  $\sX$.

\vskip 0,2 truecm
{\bf (ii)} The topological space $(\sX,\hsT)$  will be called the  {\em nonstandard 
compactification of}  $(X,T)$.

\vskip 0,2 truecm
The designation  {\em standard topology}  for  $\hsT$  arises from the fact that,  in the 
literature on nonstandard analysis,  all sets of the type  $\sG$,  where  $G\subseteq X$,  
are called  ``standard sets"  (even though  $\sG$  is, in fact, a subset of  $\sX$; 
see Definition I.2.4).                                                        

\vskip 0,2 truecm
The terminology  {\em nonstandard compactification} is justified by the following result:

\vskip 0,3 truecm\noindent
{\bf (2.5) Theorem}: Let  $(X,T)$  be a topological space and  $(\sX,\hsT)$   its nonstandard 
compactification (in the sense of the above definition).  Then:

\vskip 0,2 truecm
{\bf (i)}  Every internal subset  $A$  of  $\sX$  is compact in  $(\sX,\hsT)$.

\vskip 0,2 truecm
{\bf (ii)} $(\sX,\hsT)$  is a compact topological space  and  $(X,T)$  is a dense subspace of   
$(\sX,\hsT)$.

\vskip 0,3 truecm\noindent
{\bf Proof}: There are two ways to prove this: 1) W.A.J. Luxemburg has shown that  $\sX$  
and all internal subsets of  $\sX$  are compact with respect to the ``discrete $S$-topology" 
on $\sX$  (known also as $LS$-topology, where  $L$  stands for {\em Luxemburg)} with basic 
open sets  
$$
  ^\sigma\cP(X) = \{\sS : S\in\cP(X)\}
$$
(W.A.J. Luxemburg [14], Theorem 2.5.4, p.47 and Theorem 2.7.10, p.55). Now, the above 
statement follows from these results and the fact that the discrete $S$-topology is finer than 
the standard topology  $\hsT$.

\vskip 0,2 truecm
2) An alternative simple proof follows:  

\vskip 0,2 truecm
(i)  Let  $\{\sF_i\in\cF : i\in I\}$   be a family of basic closed sets in  $\sX$  such that the 
family  $\{\sF_i\cap\cA : i\in I\}$ has the finite intersection property. Then, by {\em Saturation  
Principle} (Chapter I, Section 2, Axiom 3),
$$
  \bigcap_{i\in I}\; \sF_i\cap\cA\not= \emptyset,   
$$
which proves that  $\cA$  is compact.

\vskip 0,2 truecm 
(ii) The compactness of  $(\sX,\hsT)$  follows from (i) as a particular case for  $\cA=\sX$. 
The original space  $(X,T)$  is a subspace of  $(\sX,\hsT)$  since $\sG\,\cap\,X = G$  for any 
$G\subseteq X$, by Corollary I.4.8, hence  $T = \{\sG\cap X:G\in T\}$. To show the denseness 
of $(X,T)$, notice that  $\sG\cap X = G\not= \emptyset$  for any basic open set  
$\sG\not=\emptyset,\; G\in T$.  The proof is complete.  $\tr$

\vskip 0,3 truecm   
It will be shown in the next chapter that  $(\sX,\hsT)$  is a  $T_0$ - space iff  $X$ is finite. 
On the other hand,  if  $X$  is finite, then  we have  $(X,T)=(\sX,\hsT)$ (Theorem I.4.6).

\vskip 0,3 truecm\noindent
{\bf (2.6) Lemma}:  For any  $H\subseteq X$  we have:

\vskip 0,2 truecm
{\bf (i)} $^\ast(\text{cl}_X H) = \text{cl}_{\sX} (\sH)$  where  ``$\text{cl}_X$" and  
``$\text{cl}_{\sX}$" are the closure operators  in  $(X,T)$  and  $(\sX,\hsT)$, respectively.  

\vskip 0,2 truecm
{\bf (ii)}  $^\ast(\text{int}_X H) = \text{int}_{\sX}(\sH)$  where ``$\text{int}_X$"  and  
``$\text{int}_{\sX}$" are the interior operators in  $(X,T)$  and  $(\sX,\hsT)$,  respectively.  

\vskip 0,3 truecm\noindent
{\bf Proof}:  We shall prove (i) only: We have
\begin{align*}
^\ast(\text{cl}_X H)\subseteq \bigcap\{\sF:\,^\ast(\text{cl}_X H)\subseteq\sF,\; X-F\in T\} = \\ 
= \bigcap\{\sF:(\text{cl}_X H)\,\subseteq F,\; X-F\in T\}  
                            =\,\bigcap\{\sF:H\subseteq F,\; X-F\in T\} =\\ 
= \bigcap\{\sF:\sH\subseteq \sF,\; X-F\in T\} = \text{cl}_{\sX} (\sH).
\end{align*}
On the other hand,  $H\subseteq\text{cl}_X H$  implies  $\sH\subseteq\,^\ast(\text{cl}_X H)$ 
which implies  $\text{cl}_{\sX} (\sH)$ $\subseteq\,^\ast(\text{cl}_X H)$  since  
$^\ast(\text{cl}_X H)$  is closed in  $(\sX,\hsT)$. The proof is complete. $\tr$

\vskip 0,3 truecm\noindent 
{\bf (2.7) Theorem} {\em (Continuity)}: Let  $(X,T)$  and  $(X',T')$  be two topological 
spaces and let $(\sX,\hsT)$  and  $(\sX',\hsT')$  be their nonstandard compactifications. 
If the function
$$
  f : (X,T)\to (X',T') \leqno{\bf (2.8)}
$$
is continuous, then  its nonstandard extension:
$$
  \ssf : (\sX,\hsT)\to(\sX',\hsT') \leqno{\bf (2.9)}
$$
is  also continuous.

\vskip 0,3 truecm\noindent 
{\bf Proof}:  For any  $G'\in T'$  we have  $\sG'\in\,^\sigma T'$  and  
$\ssf^{-1}[\sG'] =\,^\ast(f^{-1}[G'])\in\,^\sigma T$, by (Theorem I.4.12). Now, the result follows 
since  $^\sigma T$  and  $^\sigma T'$  are bases for  $\hsT$  and  $\hsT'$,  respectively.      
$\tr$

\vskip 0,3 truecm\noindent 
{\bf Note}:  It is clear that  $^\ast(f \circ g) =  \ssf \circ\,^\ast g$  and  
$^\ast(1_X) = 1_{\sX}$, so that the correspondence described above is functorial. 

\vskip 0,3 truecm\noindent 
{\bf (2.10) Theorem} {\em (Standard Part)}:  Let  $\tau$  be  the usual topology of  $\R$ 
and  $(\sR,\,^s\tau)$  be the corresponding nonstandard compactification of  $(\R,\tau)$.  
Then the standard mapping:
$$
  \text{st} : (\cF(\sR),\,^s\tau)\to (\R,\tau) \leqno{\bf (2.11)}
$$
is continuous, where  $\cF(\sR)$  denotes, as usual, the set of the finite numbers in  
$\sR$  (Definition 1.10).

\vskip 0,3 truecm\noindent 
{\bf Proof}: Let  $\alpha\in\cF(\sR)$  and  $\text{st}(\alpha) = x\in\R$.  Let  $G_x\in\tau$ 
be an open neighbourhood of  $x$  in  $(\R,\tau)$  and let  $G\in\tau$  be an open bounded 
neighourhood of  $x$  in  $(\R,\tau)$  such that  $\text{cl}(G)\subset G_x$  where  
$\text{cl}(G)$  is the closure of  $G$  in  $(\R,\tau)$.  Then  $\sG$  will be an open 
neighbourhood of  $\alpha$  in  $(\sR,\,^s\tau)$  since  $\alpha\in\mu(x)\subset\sG$.  
Moreover, we have  $\text{st}[\sG] = \text{cl}(G)\subset G_x$  by Corollary (1.12). 
That is, ``st"  is continuous at  $\alpha$  and therefore on the whole  $\cF(\sR)$.
$\tr$

\vskip 0,2 truecm
The next result is a generalization of  Theorem (2.10).

\vskip 0,3 truecm\noindent 
{\bf (2.12) Theorem}:  Let  $(X,T)$  be a regular Hausdorff space and  
$(\sX, \hsT)$  be its nonstandard compactification. Then the standard part mapping:
$$
  \text{st}:(\bigcup_{x\in X} \mu(x),\hsT)\to(X,T) \leqno{\bf (2.13)}
$$
is continuous (Definition 1.10).

\vskip 0,3 truecm\noindent 
{\bf (2.14) Notation}:  By  $C(X,\R)$  and  $C_b(X,\R)$   we shall denote the class of all  
``continuous"  and  ``continuous and bounded"  functions of the type   $f : (X,T)\to (\R,\tau)$, 
respectively, where  $(X,T)$  is a topological space and  $\tau$  is the usual topology on  $\R$.    

\vskip 0,3 truecm\noindent 
{\bf (2.15) Theorem}:  Let  $(X,T)$  be a topological space and  $(\sX,\,\hsT)$  its nonstandard 
compactification. Then, for any  $f\in C(X,\R)$  both mappings:
$$
  \ssf : (\sX,\,\hsT)\to (\sR,\,^s\tau) \leqno{\bf (2.16)}
$$
and
$$
  ^\circ f  =  \text{st}\circ\ssf :  (\widetilde X_f,\,^sT)\to(\R,\tau) \leqno{\bf (2.17)}
$$
are continuous,  where  
$$
  \widetilde X_f = \{\alpha\in\sX : \ssf (\alpha) \quad \text{is a finite number in}\,\sR\}.
    \leqno{\bf (2.18)}
$$
Also, $^\circ f$  is the unique real-valued continuous extension of  $f$  to  $\widetilde X_f$. 

\vskip 0,3 truecm\noindent
{\bf Proof}:  The continuity of  $\ssf$  follows directly from Theorem (2.7) for  $X' = \R$  and  
$T'=\tau$  and continuity of  $\;\text{st}\circ\ssf$   follows from Theorem (2.10). 
The function  $^\circ f$  is unique, since  $X$  is dense in   $\widetilde X_f$,  by Theorem (2.5). 
$\tr$

\vskip 0,3 truecm\noindent
{\bf Note}: The above result remains also true if the target space  $(\R,\tau)$ is replaced by 
a regular Hausdorff space  $(X',T')$.

\vskip 0,2 truecm
According to the notations introduced in (2.14),  $C(\sX,\R)$  will be the class of all real 
valued continuous functions defined on  $\sX$.  If  $f\in C(\sX,\R)$,  we shall denote by  
$r(f)$  the restriction of  $f$  on  $X$. 

\vskip 0,2 truecm
As a consequence of Proposition (2.5),  the next result shows that  $C(\sX,\R)$  and  
$C_b(X,\R)$  are isomorphic as rings under the restriction  map  $r$.

\vskip 0,3 truecm\noindent
{\bf (2.19) Theorem}: Let  $(X,T)$  be a topological space and  $(\sX,\,^sT)$ its nonstandard 
compactification. Then

\vskip 0,2 truecm
{\bf (i)}  For any   $f\in C(\sX,\R)$,  we have   $f = \text{st}\,\circ\,^\ast(r(f))$.

\vskip 0,2 truecm
{\bf (ii)}  $r : C(\sX,\R)\to C_b(X,\R)$   is a ring isomorphism.

\vskip 0,3 truecm\noindent
{\bf Proof}:  (i)  Given  a continuous   $f:(\sX,\,^sT)\to(\R,\tau)$,   it follows that  
$r(f)$  is continuous and bounded, since  $f$  is necessarily bounded as  $\sX$  is compact. 
Then  $f$  and  $\text{st}\circ\,^\ast(r(f))$  are two continuous functions to a Hausdorff 
space which coincide on the dense subset  $X$; hence the functions are equal on  $\sX$.   

\vskip 0,2 truecm
(ii) It is clear that  $r$  is a ring homomorphism and so is  
$s = \text{st}\,\circ\,\ast\,:\,C_b(X,\R)\to C(\sX,\R)$.  Also, (i) shows that
$s\circ r = 1$  and it is clear that   $r\circ s = 1$.  $\tr$

\vskip 0,2 truecm
Using the ``standard" theorem of M. H. Stone that  $C(X,\R)$  determines completely the 
compact Hausdorff space  $X$,  we obtain the nonstandard version which could also have been 
proved directly (with no improvements or simplifications).

\vskip 0,3 truecm\noindent
{\bf (2.20) Theorem}:  Let  $X$  and  $Y$  be compact Hausdorff spaces for which  $C(\sX,\R)$  
and  $C(^\ast Y,\R)$   are isomorphic.  Then   $X$   and  $Y$  are homeomorphic.

\vskip 1 truecm\noindent
{\bf 3. Nonstandard Hulls}

\vskip 0,2 truecm
As we show in the first section, every topological  space  $(X,T)$   can be embedded as a 
dense subspace of its nonstandard compactification   $(\sX,\,^sT)$,  having the property that any 
real valued continuous function  $f:(X,T)\to(\R,\tau)$   has a unique continuous extension:  
$$
  ^\circ f : (\widetilde X_f,\,^sT)\to(\R,\tau). \leqno{\bf (3.1)}
$$
given by  $^\circ f = \text{st}\,\circ\,\ssf$, where  $\widetilde X_f$  is defined in (2.18). 
The space  $(\sX,\,^sT)$  is  Hausdorff only when   $X$  is finite and Hausdorff.

\vskip 0,2 truecm
Following the ``nonstandard hull construction" (W.A.J. Luxemburg [14]),  we shall 
consider  the factor space :
$$
   \hat X_\Phi =  \widetilde X_\Phi\,/\sim_\Phi \leqno{\bf (3.2)}
$$
indentifying points of a given subset  $\widetilde X_\Phi$  of  $\sX$  under an equivalence 
relation  ``$\sim_\Phi$". We shall specify   $\widetilde X_\Phi$   and  $\sim_\Phi$   in terms 
of a given family of real valued continuous functions  $\Phi\subseteq C(X,\R)$.

\vskip 0,3 truecm\noindent
{\bf (3.3) Definition}:  Let  $\Phi\subseteq C(X,\R)$. Then :

\vskip 0,2 truecm
{\bf (i)}  $\quad \widetilde X_\Phi$  consists of all points  $\alpha$  in  $\sX$  such that  
$\ssf(\alpha)$  is a finite number in  $\sR$  for all  $f\in\Phi$. The points in  
$\widetilde X_\Phi$  will be called ``$\Phi$-finite".

\vskip 0,2 truecm
{\bf (ii)}   Two points  $\alpha$  and  $\beta$  in  $\tX_\Phi$,  are called $\Phi$-equivalent, 
written as  $\alpha\sim_\Phi \beta$,  if   $\,\ssf(\alpha)\approx\ssf(\beta)\,$  for  all  $\,f\in\Phi\,$,  
where  $\;\approx\;$  is the infinitesimal relation in  $\sR$. 

\vskip 0,2 truecm
{\bf (iii)} The factor space  $\tX_\Phi$  will be given the {\em quotient topology}  $\hT$. 
The corresponding topological space  $\,(\hX_\Phi,\hT)\,$  will be called$\;$ ``the nonstandard 
$\;\Phi$-hull of $(X,T)$".

\vskip 0,2 truecm
{\bf (iv)} For every  $f\in\Phi$,  there is a well defined  mapping 
$$
  \hf : \hX_\Phi\to\R,
$$ 
given by  $\hf\circ q =\,^\circ f$,  where  $q$  is the quotient mapping from  $\tX_\Phi$ 
onto  $\hX_\Phi$.  

\vskip 0,2 truecm
The following result establishes a connection between the monads of the space  $(X,T)$  
and the equivalence relation  $\sim_\Phi$.    

\vskip 0,3 truecm\noindent
{\bf (3.4) Lemma}:   $\mu(x)\subseteq q(x)$  for any  $x\in X$.  When  the family  $\Phi$  
distinguishes  points and closed sets in  $X$, then  $\mu(x) = q(x)$  for  all  $x$  in  $X$.

\vskip 0,3 truecm\noindent
{\bf Proof}:  $\mu(x)\subseteq q(x)$  follows immediately from the fact that all  $f$  in  $\Phi$  
are continuous and therefore,  $\ssf(\alpha)\approx\ssf(x)= f(x)$  for all  $\alpha\in\mu(x)$. 
Let  $\Phi$  distinguish  points and  closed sets,  i.e. for each closed  $F\subset X$  and  
$x\in X-F$,  $g(x)\notin\text{cl}\,g[F]$  for some  $g\in\Phi$. Let  $\alpha\in q(x)$,  i.e.  
$\alpha\sim_\Phi x$,  which means   $\ssf(\alpha)\approx f(x)$  for  all  $f\in\Phi$. 
We have to show that  $\alpha\in\mu(x)$. Suppose (for contradiction) that  
$\alpha\not\in\sG$  
for some open neighbourhood  $G$  of  $x$  in  $(X,T)$  and choose   $F = X-G$. 
There exists   $g\in\Phi$ which distinguishes  $x$  from  $X-G$  in the sense that 
$g(x)\notin\text{cl}_{\Bbb R}\, (g[X-G]$). On the other hand, we have 
$$
  \text{cl}_{\Bbb R}\;(g[X-G]) = \{y\in\R : y\approx\,^\ast g(\beta)\quad\text{for some}\quad
  \beta\in\sX-\sG\},
$$
by Corollary 1.12, since $^\ast(g[X - G]) =\, ^\ast g[\sX - \sG]$, by Theorem I.4.12. 
It follows  $^\ast g(\alpha)\not\approx g(x)$,  contradicting   $\alpha\sim_\Phi x$.
$\tr$

\vskip 0,5 truecm
It should be noted that not all topological spaces  $(X,T)$  admit  families of continuous 
real valued functions  $\Phi$  which distinguish  points and closed sets. The spaces which 
admit  $\Phi$  with this property are the completely regular  ones (J.L. Kelley [12]).

\vskip 0,3 truecm\noindent
{\bf (3.5) Theorem}:  The quotient mapping  $q:\tX_\Phi\to(\hX_\Phi,\hT)$  maps  $X$  onto 
a dense subset of  $(\hX_\Phi,\hT)$.

\vskip 0,3 truecm\noindent
{\bf Proof}:  $X$  is dense in  $(\sX,\,^sT)$,  by Theorem (2.5), hence, dense in  $\tX_\Phi$.
Therefore, $q[X]$ is dense in  $q[\tX_\Phi]=\hX_\Phi$, by continuity. $\tr$

\vskip 0,3 truecm\noindent 
{\bf (3.6) Theorem}: For every  $f\in \Phi$,  the mapping
$$
  \hf : (\hX_\Phi,\hT)\to (\R,\tau) \leqno{\bf (3.7)}
$$
is continuous and  $\hf$ is the unique real-valued continuous extension of $f$ to $\hX_\Phi$,
in the sense that
$$
   f(x) =  \hf(q(x)),\qquad  x\in X. \leqno{\bf (3.8)}
$$

\vskip 0,3 truecm\noindent
{\bf Proof}: As remarked above  $\hf$  is well defined on   $\hX_\Phi$.  Since  $(\hX_\Phi,\hT)$ 
has the quotient topology induced by  $q,\hf$   is continuous iff  $\hf\circ q$  is continuous. 
Now,  $\hf\circ q =\,^\circ f$  is continuous by Theorem (2.15). Finally,  
$f(q(x))=\,^\circ f(x)=(\text{st}\circ\ssf)(x) = \text{st}(\ssf(x)) = f(x)$. 
The function  $\hf$  is unique, since $q[X]$  is a dense subset of  $\hX_\Phi$, by Theorem 
(3.5). The proof is complete. $\tr$\

\vskip 0,3 truecm\noindent
{\bf (3.9) Theorem}:  Let  $X,T)$  be a topological space and  $\Phi\subseteq C(X,\R)$. 
Then the corresponding $\Phi$-hull  $(\hX_\Phi,\hT)$  is a Hausdorff space.

\vskip 0,3 truecm\noindent
{\bf Proof}:  Let  $a,b\in\hX_\Phi$  be two distinct points.  Then there are points  
$\alpha,\beta$  in  $\tX_\Phi$  such that  $q(\alpha) = a$,  $q(\beta)=b$  and a function 
$f\in\Phi$  for which   $\ssf(\alpha)\approx\ssf(\beta)$.  Then  $\hf(a)\not=\hf(b)$  in  $\R$, 
so there are disjoint open sets in  $\R$,  $U$  and  $V$, with $\hf(a)\in U$,  $\hf(b)\in V$.  
Now  $\hf$  is continuous, so  $a\in\hf^{-1}[U]$,   $b\in\hf^{-1}[V]$, as required. $\tr$

\vskip 0,3 truecm
We consider now some particular cases for the family  $\Phi$  and the initial topological 
space  $(X,T)$. First, we obtain the Hausdorff compactifications of  $(X,T)$.

\vskip 0,3 truecm\noindent
{\bf (3.10) Theorem}: If  $\Phi\subseteq C_b(X,\R)$, then:

\vskip 0,2 truecm
{\bf (i)}  $\,\tX_\Phi = \sX$.

\vskip 0,2 truecm
{\bf (ii)} $(\hX_\Phi,\hT)$  is a compact space containing a continuous image of  $(X,T)$.  

\vskip 0,3 truecm\noindent
{\bf Proof}: (i) $\Phi\subseteq C_b(X,\R)$  implies  $\tX_\Phi = \sX$  since for bounded functions  $f$
all values of  $\ssf$  are finite numbers in  $\sR$. 

\vskip 0,2 truecm
(ii) follows immediately from Theorem (3.5) and the fact that the continuous image  
$q[\sX]$  of a compact space  $\sX$  is compact.   $\tr$

\vskip 0,3 truecm\noindent
{\bf (3.11) Corollary}: Let $(X,T)$  be a completely regular Hausdorff space and let the family 
$\Phi$ distinguish the points and the closed sets in  $(X,T)$.  Then  $(X,T)$  is homeomorphic 
to its image in  $(\hX_\Phi,\hT)$,  in symbols,  $X\subseteq\hX_\Phi$,  and for any  $f$ in  $\Phi$   
we have
$$
 \hf(x) = f(x),\qquad\qquad\qquad\qquad x \in X. \leqno{\bf (3.12)}
$$

\vskip 0,3 truecm\noindent
{\bf Proof}:  Since  $\Phi$  distinguishes the points and the closed sets in  $(X,T)$, we have 
$q(x) = \mu(x)$  for any   $x\in X$, by Lemma (3.4). Also  $q(x) = q(y)$  if and only if 
$x=y$  for any  $x,y\in X$, since $(X,T)$  is Hausdorff. That means that the quotient mapping  
$q$  is one to one. Let  $s$   from  $X\subset\hX_\Phi$  to  $X\subset\tX_\Phi$  be the inverse 
of  $q$. Now,  $(X,T)$  is a completely regular space so  $s$  is continuous if and only if   
$f\circ s$   is continuous  for all   $f\in\Phi$.  But  $f\circ s = f\mid X$. The formula (3.12) 
follows immediately from  (3.8).   The proof is complete.  $\tr$

\vskip 0,2 truecm
It is instructive to illustrate the above proceedure for special families  $\Phi$.

\vskip 0,5 truecm\noindent
{\bf (3.13) Examples}:       
\vskip 0,2 truecm
{\bf 1.} If  $\Phi$  is empty,  then  $\tX_\Phi = \sX$,  all points are equivalent and $\hX_\Phi$  
reduces to a single point.

\vskip 0,2 truecm
{\bf 2.}  Consider  $\Phi = \{\text{id}\}$,  where  $\text{id}:(\R,\tau)\to (\R,\tau)$   
is the identity map.  Then  $\tX_\Phi$  is   $\cF(\sR)$   and   $\alpha\sim_\Phi \beta$  iff   
$\alpha,\beta\in\mu(x)$   for some  $x\in\R$. We have  $(\hX_\Phi,\hT)=(\R,\tau)$.

\vskip 0,2 truecm 
{\bf 3.}  Again, consider  $(\R,\tau)$  and  $\Phi = \{\sin x, \cos x\}$. Then   
$\tX_\Phi=\sR$, $\alpha\sim\beta$  iff   $|\alpha-\beta|\approx 2k\pi$  for some  $k\in\Bbb Z$  
and   $q[\R]$  is, topologically, the circle $\{(x,y):x^2+y^2 = 1\}$  with the Euclidean topology. 
Thus, $\hX_\Phi = q[\R]$. 

\vskip 0,2 truecm
{\bf 4.}  If  $\Phi$  consists of all real valued bounded functions on  $X$,  then  
$(\hX_\Phi,\hT)$  is the Stone-\v Cech compactification of  $(X,T)$.

\vskip 0,2 truecm
{\bf 5.}  If  $\Phi$  consists of  all real valued continuous functions on  $X$,  then 
$(\hX_\Phi,\hT)$  is the Hewitt realcompactification of  $(X,T)$.

\vskip 0,2 truecm
Both  no. 4.  and  no. 5.  will be established in the next sections.

\vskip 0,5 truecm\noindent
{\bf (3.14) Theorem}: If  $\Phi$  is  $C_b(X,\R)$  or  $C(X,\R)$,  then  $\Phi$  and 
$C(\hX_\Phi,\R)$ are isomorphic as rings (for the notation see (2.14)).

\vskip 0,3 truecm\noindent
{\bf Proof}:  For each  $f$  in  $\Phi$  we have shown that  $\hf$   is continuous and 
$f=\hf\circ q$  on  $X$. This defines a map  $\varphi : \Phi\to C(\hX_\Phi,\R)$. This map 
is injective, since  $\hf_1 = \hf_2$   gives   $\hf_1\circ q = \hf_2\circ q$, i.e. $f_1=f_2$.  
It is surjective, for suppose  $g:(\hX_\Phi,\hT)\to (\R,\tau)$. Let  $f$  be the restriction  
of  $g\circ q$  to  $X$. We show that  $g = \hf$.  This will follow from  $g\circ q=\hf\circ q$  
on  $X$. For  $x\in X$   we have  $(g\circ q)(x)= f(x)$,  by the definition of  $f$;  
also $(\hf\circ q)(x)=f(x)$  by definition of  $\hf$.  Hence  $g = \hf$.  Finally,  $\varphi$ 
is a ring isomorphism. We verify only one property:  $\varphi(f_1+f_2)=\varphi(f_1)+\varphi(f_2)$
iff  $\varphi(f_1+f_2)\circ q = \varphi(f_1)\circ q + \varphi(f_2)\circ q$  on  $X$
iff  $(f_1+f_2)^\wedge\,\circ q = \hf_1\circ q + \hf_2\circ q$  on  $X$  iff  $f_1+f_2 = f_1+f_2$  
on  $X$. The proof is complete.  $\tr$

\vskip 1 truecm\noindent
{\bf 4. Stone - \v Cech Compactification: The Case  $\Phi = C_b(X,\R)$} 

\vskip 0,3 truecm
Let  $(X,T)$  be a topological space and let  $\Phi$, which appears  in Definition (3.3), be the 
class of continuous bounded real valued functions defined on  $X$,  i.e.  $\Phi = C_b(X,\R)$. 
In this particular case we have  $\tX_\Phi = \sX$, by Proposition (3.10). Throughout this section 
we shall write simply  $\sim$  and  $\tX$  (suppressing the index  $\Phi$)  instead of the more 
precise ``$\sim_\Phi$ and  $\hX_\Phi$  for  $\Phi = C_b(X,\R)$", respectively. In this notation, 
for the nonstandard hull we have: $\hX = \sX /\sim$,  where  $\alpha\sim\beta$  in  $\sX$  iff 
$\ssf(\alpha)\approx\ssf(\beta)$  for all  $f$ in  $C_b(X,\R)$. Let  $(\hX,\hT)$ be the 
corresponding topological space (Definition (3.3)).

\vskip 0,5 truecm\noindent
{\bf (4.1) Theorem}: $(\hX,\hT)$  coincides with the Stone-\v Cech compactification $\beta X$
of  $(X,T)$.

\vskip 0,3 truecm\noindent
{\bf Proof}: The space  $\hX$  is Hausdorff and compact, by Theorem (3.9) and  Theorem (3.10), 
respectively. Also, $q[X]$ is a dense subset of  $\hX$  by Theorem (3.5), which immediately 
implies the uniqueness of all continuous extensions  $\hf$  of  $f$ (Theorem (3.6)). These 
properties characterize   $\beta X$.  $\tr$

\vskip 0,5 truecm\noindent
{\bf (4.2) Corollary}: {\em (Completely Regular Hausdorff Space)}:  Let $(X,T)$ be a completely 
regular Hausdorff space. Then  $(X,T)$  is homeomorphic to its image in  $(\hX,\hT)$, in symbols, 
$X\subseteq\hX$, and for any  $f$  in  $C_b(X,\R)$ we have  $\hf(x) = f(x)$  for all   $x\in X$.

\vskip 0,3 truecm\noindent
{\bf Proof}: Since  $(X,T)$  is completely regular, the family  $C_b(X,\R)$  distinguishes the 
points and closed sets in  $X$. Now, the result follows directly from Corollary (3.11). $\tr$

\vskip 0,2 truecm
Compared with other nonstandard expositions of the Stone-\v Cech compacti\-fication ([4], [10], 
[14], [15], [17], [18], [21], [22]) we wish to emphasize that {\em we do not use the weak 
topology neither} on  $X$,  nor on  $\hX$. Continuous functions from  $C_b(X,T)$  are only used 
to define the equivalence relation in  $\sX$. \\

{\bf 5. All Compactifications}

\vskip 0,2 truecm
Let  $(X,T)$  be a topological space.  A compact Hausdorff space  $(K,L)$  is a 
``compactification of  $(X,T)$"  if there is a continuous function  $\psi : (X,T)\to (K,L)$
such that   $\psi[X]$  is dense in  $(K,L)$.

\vskip 0,2 truecm
This definition includes the more familiar and restrictive definition of a Hausdorff 
compactification of a completely regular space  $(X,T)$  as one that contains  $(X,T)$  
as a dense subspace.

\vskip 0,2 truecm
The purpose of this section is to show that all Hausdorff compactifications of  $(X,T)$  
can be obtained as nonstandard hulls in the manner described in Section 3. 

\vskip 0,2 truecm
The question of obtaining all compactification of a given completely regular Hausdorff 
space, in the more restricted sense mentioned above, has been considered by K.D. Stroyan 
[21] in terms of an infinitesimal relation induced in the category of totally bounded uniform 
spaces. In our approach the relation is purely topological and the given compactification is  
$\sX/\sim_\Phi$  for suitable  $\Phi$.

\vskip 0,2 truecm
Consider a Hausdorff compactification  $(K,L)$  of  $(X,T)$  with a continuous map  
$\psi : X\to K$  with dense range, we shall keep  $X,K$  and  $\psi$  fixed throughout 
the following discussion. 

\vskip 0,2 truecm
There is the continuous extension  $^\ast\psi : (\sX,\,^sT)\to (\sK,\,^sL)$  (Proposition (1.7)) 
and the continuous standard map function  $\text{st}_K : (\sK,\,^sL)\to (K,L)$  
(Definition (1.10)),  so that  $\Psi:(\sX,\,^sT)\to (K,L),\Psi = \text{st}_K \circ\,^\ast\psi$,  
gives a continuous extension of  $\psi$ (on $X$) to $\sX$. Moreover, if  $f:(K,L)\to(\R,\tau)$  
is continuous function, then  $\ssf:(\sK,\,^sL)\to(\sR,\,^s\tau)$  is continuous and   
$f\circ\text{st}_K = \text{st}_\Bbb R\circ\ssf$, since  $\ssf[\mu(x)]\subseteq m(f(x))$, 
where  $\mu$  and  $m$  are the monads of  the spaces  $(X,T)$  and  $(\R,\tau)$  respectively.  
This situation is best summarized in the commutative diagram that follows:

$$\begin{array}{lllll}
   & \phantom{a}^\ast\psi & & \phantom{a}^\ast f &  \\
   \phantom{a}^\ast  X
   {\mbox{\setlength{\unitlength}{1cm}\put(0.3,-0.1){\vector(1,-1){1.3}}}}
  & \overrightarrow{\phantom{aaaaa}} &\phantom{}^\ast K 
& \overrightarrow{\phantom{aaaaa}} &\phantom{}^\ast \R\\[-1.2cm]
\phantom{a}^\ast \Bigg\uparrow & \ \ \ \ \ \ \Psi  
 & \ \,\Bigg\downarrow & \hskip-0,5cm \text{st}_K   & \ \,\Bigg\downarrow \text{st}_{\Bbb R}\\[0.5cm]
  \phantom{aa} X & \overrightarrow{\phantom{aaaaa}}  & \ \,K 
&\overrightarrow{\phantom{aaaaa}} & \ \,\R\\[-0.2cm]
  & \ \ \psi & & \ \ f & 
\end{array} \leqno{\bf (5.1)}$$

\noindent						
{\bf (5.2) Definition} {\em (The Family  $\Phi$)}: Let  $\Phi$  consist of all  
$f\circ\psi,\; f\in C(K,\R)$  (for the notation see (2.14)).

\vskip 0,1 truecm
Thus, $\Phi$  consists of all real valued continuous  $g$  on  $(X,T)$  which have  an 
``extension"  $f$  to  $(K,L)$  in the sense  $f$  is continuous and  $g = f \circ\psi$.

\vskip 0,1 truecm
Observe that  $\Phi\subseteq C_b(X,\R)$, so that  $\Phi$  determines an equivalence relation  
``$\sim_\Phi$"  and  $\tX_\Phi = \sX$  such that  $(\hX,\hT,)$  is a Hausdorff  compactification 
of  $(X,T)$  where  $(X,T)\to(\hX,\hT)$  is given by the restriction of  
$q:(\sX,\,^sT)\to(\hX,\hT)$  on  $X$  (Theorem (3.9) and Theorem (3.10)). We show that  
$(\hX,\hT)$  is homeomorphic to  $(K,L)$.

\vskip 0,5 truecm\noindent	
{\bf (5.3) Lemma}: For  $f:(K,L)\to(\R,\tau)$  and  $\Psi:(\sX,\,^sT)\to(K,L)$,  as above, 
we have :                                                  
$$
  f\circ\Psi  =  \text{st}_\Bbb R\,^\ast(f\circ\psi). \leqno{\bf (5.4)}
$$

\vskip 0,3 truecm\noindent	
{\bf Proof} :  $f\circ\Psi = f\circ\text{st}_K\circ\,^\ast\psi = 
\text{st}_\Bbb R\circ\ssf\circ\,^\ast\psi = \text{st}_\Bbb R\,^\ast(f\circ\psi)$.  $\tr$

\vskip 0,5 truecm\noindent	
{\bf (5.5) Lemma}:  $\alpha\sim_\Phi\beta$   if and only if   $\Psi(\alpha) = \Psi(\beta)$.

\vskip 0,3 truecm\noindent	
{bf Proof}:  Suppose  $\Psi(\alpha)\not= \Psi(\beta)$.  Since  $(K,L)$  is compact Hausdorff, 
there is  
$$
  f:(K,L)\to([0,1],\tau)
$$  
such that  $f(\Psi(\alpha))= 0, f(\Psi(\beta))= 1$. But then  
$\text{st}_\Bbb R\,^\ast(f\circ\psi)(\alpha)\not= \text{st}_\Bbb R\,^\ast(f\circ\psi)(\beta)$,  
which contradicts  $\alpha\sim_\Phi\beta$.  The converse is clear.  $\tr$

\vskip 0,5 truecm
The proposition above shows that there is a well defined map  $\chi:\hX\to K$  given by  
$\chi\circ q = \Psi$.  Since  $\hT$  is the quotient  topology induced by  $q$, we have:

\vskip 0,5 truecm\noindent	
{\bf (5.6) Proposition}: There is a continuous map  $\chi : (\hX,\hT)\to(K,L)$  such that  
$\chi\circ q = \Psi$.

\vskip 0,3 truecm\noindent	
{\bf Note}:  The mapping $\chi$  obtained above satisfies  $\chi\circ q = \psi$  on  $X$  
so it must be the Stone extension  $\psi^\beta$   of  $\psi : X\to K$  
(see L. Gillman and  M. Jerison  [3]), by uniqueness of that extension.

\vskip 0,5 truecm\noindent	
{\bf (5.7) Theorem}:  $(\hX,\hT)$  and  $(K,L)$  are homeomorphic.

\vskip 0,5 truecm\noindent	
{\bf Proof}:  Clearly,  $\chi(q[X]) = \psi[X]$  is dense in  $(K,L)$, so that  $\chi$  is 
surjective. Since   $(\hX,\hT)$  and  $(K,L)$  are compact Hausdorff, it only remains to 
show that $\chi$  is injective. But this follows from the fact that  $\alpha\sim_\Phi\beta$  
iff  $\Psi(\alpha) = \Psi(\beta)$.    $\tr$      

\vskip 0,5 truecm
The mapping  $\Psi:\sX\to K$  allows a simple description for zero sets  $Z(f)$  which we 
shall give in what follows. In particular, when  $X$  is completely regular and Hausdorff 
and  $K$  is the Stone-\v Cech compactification, we obtain a description of  
$\text{cl}_{\beta X}\,Z(f)$  and  $Z(f^\beta)$  for  $\in C_b(X,\R)$. As above, we assume that 
$(K,L)$  is a Hausdorff compactification of  $(X,T)$  with  $\psi : (X,T)\to(K,L)$, and  
$\psi[X]$  is dense in  $K$. 

\vskip 0,5 truecm\noindent	
{\bf (5.8) Propositon}:   Let  $g : (X,T)\to(\R,\tau)$  be  such  that  there  is  an  extension
$f:(K,L)\to(\R,\tau)$  with  $g = f\circ\psi$.  Then
$$
  \Psi^{-1}[Z(f)] = \{\alpha\in\sX : \ssf(\alpha)\approx 0\}. \leqno{\bf (5.9)}
$$

\vskip 0,2 truecm\noindent	
{\bf Proof}: $^\ast g(\alpha)\approx 0\Leftrightarrow (\ssf\circ\,^\ast\psi)(\alpha)\approx
0\Leftrightarrow\text{st}_\Bbb R\,(\ssf\circ\,^\ast\psi)(\alpha) = 0 \Leftrightarrow
(f\circ\Psi)(\alpha) = 0 \Leftrightarrow \Psi(\alpha)\in Z(f)$.    $\tr$

\vskip 0,5 truecm
When  $(X,T)$   is completely regular and Hausdorff and  $K$  is  $\beta X$, we regard $\psi$
as the identity on  $X$  and  so,  $\Psi$  is  $q: \sX\to X$  and the statement above is:   
$$
  Z(g^\beta) =  q [\{\alpha\in\sX :\,^\ast g(\alpha)\approx 0\}]. \leqno{\bf (5.10)}
$$
It is of interest to observe that it is pointed out in  $L$. Gillman and M. Jerison's book [3]  
that  $Z(g^\beta)$  need not be of the form  $\text{cl}_{\beta Z}(h)$,
$h\in C_b(X,T)$, and  that  $Z(g^\beta)$  is always a countable intersection of sets of the form  
$\text{cl}_{\beta X}\,Z(f)$ ([6],  6E and, also,  8D). The formula above gives the precise 
description of  $Z(g^\beta)$. When  $g:\N\to\R$  is  $g(n) = 1/n$,  then  $Z(g^\beta)$  is  
the image of all infinitely large natural numbers under  $q$.

\vskip 0,5 truecm\noindent
{\bf (5.11) Proposition}: Let  $g:(X,T)\to(\R,\tau)$  be a bounded function. Then
$$
  \text{cl}_K(\psi[Z(g)]) = \Psi[\sZ(g)]. \leqno{\bf (5.12)}
$$

\vskip 0,3 truecm\noindent
{\bf Proof}:  It is clear that  $\psi[Z(g)]\subseteq\Psi[Z(g)]\subseteq\Psi[\sZ(g)]$. 
The last set is compact, hence closed in  $K$,  so that  
$\text{cl}_K\,(\psi[Z(g)])\subseteq\Psi[\sZ(g)]$.  Conversely,  $\sZ(g)=\text{cl}_{\sX}(\sZ(g))$,
by Lemma (2.6),  so that  $\Psi[\sZ(g)] = \Psi[\text{cl}_{\sX}(\sZ(g))]\subset
\text{cl}_K\,[\Psi(Z(g))] = \text{cl}_K\,(\psi[Z(g)])$.     $\tr$

\vskip 0,5 truecm
As before, when $(X,T)$ is completely regular and  Hausdorff  and  $K$  is $\beta X$, we have:
$$
  \text{cl}_{\beta X}\,(Z(g)) = q [\sZ(g)] = q [\{\alpha\in\sX :\,^\ast g(\alpha) = 0\}].
    \leqno{\bf (5.13)}
$$

This formula, combined with a classical standard characterization of  $\beta X$  (L. Gillman and M. Jerison [3], (6.5), IV)  gives a nonstandard characterization of  $\beta X$  which we 
formulate for completely regular spaces.

\noindent
{\bf (5.14) Proposion}:  Let  $(X,T)$  be a completely regular Hausdorff space and  $(K,L)$  
a compact Hausdorff space containing  $(X,T)$  as a dense subspace.  Then  $(K,L)$  is the  
Stone-\v Cech compactification of  $(X,T)$  if and only if  for any zero sets  $Z_1$, $Z_2$  
in  $X$  we have:
$$
  \Psi[\sZ_1\cap\sZ_2]  =  \Psi[\sZ_1]\cap\Psi[\sZ_2]. \leqno{\bf (5.15)}
$$

\vskip 0,2 truecm\noindent
{\bf Proof}: If  $(K,L)$  is  $\beta X$,  then for zero sets  $Z_1$, $Z_2$  we have  
$$
  \text{cl}_{\beta X}\,(Z_1\cap Z_2) = \text{cl}_{\beta X}\,Z_1\cap\text{cl}_{\beta X}\,Z_2
$$
(L. Gillman and M. Jerison [3], (6.5), Compactification Theorem).  Hence  
$$
  \Psi[^\ast(Z_1\cap Z_2)]  =  \Psi[\sZ_1]\cap \Psi[\sZ_2].
$$
Now   $^\ast(Z_1\cap Z_2) = \sZ_1\cap\sZ_2$  and the result follows. The proof of the converse 
is similar.   $\tr$\\

\noindent
{\bf 6. Hewitt Realcompactification : The Case $\Phi = C(X,\R)$}

\vskip 0,2 truecm\noindent 
Let  $(X,T)$  be a topological space.  We mentioned in Example (3.13) - no.5.,  that if we 
put  $\Phi= C(X,\R)$  in Definition (3.3), the corresponding  $\Phi$-hull will coincide 
with the Hewitt realcompactification of  $(X,T)$  (L. Gillman and M. Jerison [3]). 
We now discuss this important case in detail. 

\vskip 0,2 truecm\noindent
We shall write simply  $\tX$, $\sim$, and  $\hX$ (suppressing the index $\Phi$)  instead of 
the more precise ``$\tX_\Phi,\,\sim_\Phi$,  and  $\hX_\Phi$,  for  $\Phi= C(X,\R)$", respectively, 
throughout the following discussion. ({\em Warning}:  $\tX$  should not be confused with the set 
of the nearstandard points of $\sX$). Since we have to extend all (not only the bounded) 
continuous functions to the new space, we have to select for the set of the $\Phi$-finite 
points some proper subset of  $\sX : \tX$  is the set of all points  $\alpha$  in  $\sX$  
for which  $\ssf(\alpha)$  is a finite number in $\sR$  for all  $f$  in  $C(X,\R)$. For the 
nonstandard hull we have  $\hX = \tX/\!\!\sim\,,$  where  $\alpha\sim\beta$  in  $\tX$  if and only if 
$\ssf(\alpha)\approx\ssf(\beta)$  for all  $f$  in  $C(X,\R)$.  Let $(\hX,\hT)$ be the 
corresponding topological space (Definition (3.3)).

\vskip 0,2 truecm
Recall that a topological space  $(X,T)$  is called  ``realcompact" if for every nontrivial 
ring homomorphism  $\pi : C(X,\R)\to\R$  there is  $x\in X$  such that ``$\pi(f) = 0$  iff  
$f(x) = 0$" for all  $f\in C(X,\R)$ or, equivalently, if ``every real maximal ideal of 
$C(X,\R)$ is fixed" (L. Gillman and M. Jerison [3]). 

\vskip 0,3 truecm\noindent
{\bf (6.1) Lemma}: Let  $\pi: C(X,\R)\to\R$   be a nontrivial ring homomorphism. Then, there 
exists  $\alpha$  in   $\tX$   such that  $\pi(f) = \ssf(\alpha)$   for all  $f\in C(X,\R)$.

\vskip 0,3 truecm\noindent
{\bf Proof}: The family of internal subsets of  $\sX : A_f = \ssf^{-1}[\{0\}]$,  
$f\in\text{ker}\,\pi$, 
has the finite intersection property. Indeed,  $f^{-1}[\{0\}]\subseteq\ssf^{-1}[\{0\}]$   and,  
on the other hand, $f^{-1}[\{0\}] = \emptyset$  implies that  $f$  is invertible in  $C(X,\R)$  
which contradicts  $f\in\text{ker}\,\pi$. So that,  $\ssf^{-1}[\{0\}]\not=\emptyset$  and moreover, 
we have:
$$
  \ssf^{-1}[\{0\}]\cap\,^\ast g^{-1}[\{0\}]=\,^\ast(f^{-1}[\{0\}]\cap g^{-1}[\{0\}])\supseteq
  \,^\ast(f^2 + g^2)^{-1}[\{0\}]\not= \emptyset.
$$

\vskip 0,2 truecm\noindent
By the {\em Saturation Principle}  (Chapter I, Section 2, Axiom 3), there exists  $\alpha\in\sX$  
such that $\ssf(\alpha) = 0$  for all  $f\in\text{ker}\,\pi$. Taking into account also that  
$\text{ker}\,\pi$  is a maximal ideal of  $C(X,\R)$, we get                   
$$
  \text{ker}\,\pi = \{f\in C(X,\R)\mid \ssf(\alpha) = 0\}.
$$

\vskip 0,1 truecm\noindent
Now, for  $f\in C(X,\R)$,  we have   $\pi(f) = c\in\R$. Then, we have $f-c\in\text{ker}\,\pi$  
so,  $\ssf(\alpha) = c = \pi(f)$. Since  $c$  is  a real number,  $\alpha\in\tX$. 
The proof is complete.   $\tr$

\vskip 0,3 truecm\noindent
{\bf Note}:  The result of the above lemma is related to results in  (J.C. Dyre [2], Theorem (3.3)).   
The difference with  Dyre's work consists in our restriction to real maximal ideals of  $C(X,\R)$ 
only and, hence, the localization of  $\alpha$  in  $\tX$  which is essential for our discussion.

\vskip 0,3 truecm\noindent
{\bf (6.2) Theorem}:  $(\hX,\hT)$  is realcompact.

\vskip 0,3 truecm\noindent
{\bf Proof}: Let  $\pi : C(\hX,\R)\to\R$   be a nontrivial ring homomorphism. Then, define 
$\varphi : C(X,\R)\to C(\hX,\R)$  by   $\varphi(f) = \hf$  (Definition (3.3))   and observe 
that the map:  $\pi\circ\varphi : C(X,\R)\to\R$  is also a nontrivial ring homomorphism. 
Then, by Lemma (6.1), there is $\alpha\in\tX$  such  that  $(\pi\circ\varphi)(f)=\ssf(\alpha)$  
for all  $f\in C(X,\R)$  which means  $\pi(\hf) = \ssf(\alpha) = \hf(q(\alpha))$  for all 
$f\in C(X,\R)$. Taking into acount  Theorem (3.14), we get that $\pi(\hf)= 0$ iff 
$\hf(\alpha)=0$  for all  $\hf\in C(\hX,\R)$ where  $a = q(\alpha)$. The proof is complete.  $\tr$ 

\vskip 0,3 truecm\noindent
{\bf (6.3) Lemma}: Let  $f\in C(X,\R)$  and  $\alpha$  and  $\beta$  in  $\sX$  be such that 
$\ssf(\alpha)\approx\ssf(\beta)$. If  $\alpha\in\tX$, then there is a continuous function  
$g : X\in[0,1]$ such that $\sg(\alpha) = 0$  and   $\sg(\beta) = 1$.

\vskip 0,3 truecm\noindent
{\bf Proof}: Since  $\alpha\in\tX$, the value  $\ssf(\alpha)$  is a finite number in  $\R$  so, 
whether  $\ssf(\beta)$  is infinitely large or not, there are open sets  $U$, $V$ in $\R$  whose 
closures are disjoint and  $\ssf(\alpha)\in\sU$  and  $\ssf(\beta)\in\sV$.  Let 
$\varphi:\R\to[0,1]$  be continuous and such that $\varphi$  is  $0$  on  $U$  and  $1$ on  $V$. 
The function  $g = \varphi\circ f$  has the required properties, since:
$$
  ^\ast\varphi^{-1}[\{0\}] = \,^\ast(\varphi^{-1}[\{0\}])\supseteq\sU\ni\ssf(\alpha),
$$
i.e. $\sg(\alpha) = \,^\ast\varphi(\ssf(\alpha)) = 0$  and, similarly,  
$\sg(\beta) = \,^\ast\varphi(\ssf(\beta)) = 1$.   $\tr$

\vskip 0,2 truecm
Somewhat surprisingly, it is possible to prove that  $(\hX,\hT)$ is completely regular. 
The argument uses the compactness of  $(\sX,\,^sT)$  (Theorem (2.5)).

\vskip 0,3 truecm\noindent
{\bf (6.4) Proposition}:  $(\hX,\hT)$   is a completely regular space.

\vskip 0,3 truecm\noindent
{\bf Proof}: Let  $a\in\hX$  and let  $F\subseteq\hX$  be a closed set not containing $a$. Then 
$q^{-1}[F]$  is a closed subset of  $\tX$, so there is a closed set  $K$ in  $\sX$  such that 
$K\cap\tX = q^{-1}[F]$. Since  $\sX$  is compact, $K$  is also compact in  $\sX$. Then, let 
$\alpha\in\tX$  be such that  $q(\alpha) = a$. Clearly  $\alpha\notin K$. Moreover, for each 
$\beta\in K$  there exists  $f_\beta \in C(X,\R)$  such that  
$\ssf_\beta(\alpha)\not\approx\ssf_\beta(\beta)$. For suppose not, then we obtain  $\beta\in\tX$   
and  $\alpha\sim\beta$,  i.e.  $a = q(\alpha)\in F$, a contradiction. By Lemma (6.3), we may 
assume that $0\le f_\beta\le 1$  and  $\ssf_\beta(\alpha) = 0$, $\ssf_\beta(\beta) = 1$. 
Then the sets $(\ssf_\beta)^{-1}[^\ast(3/4,1]]$ cover  $K$, so there are finitely many such sets  
$(\ssf_r)^{-1}[^\ast(3/4,1]]$, $r = 1,2,\dots,n$, which cover  $K$.  Also, $\ssf_r(\alpha) = 0$,  
$r = 1,2,\dots,n$. Let $g = \sup\{f_r : 1\le r\le n\}$. Then,  
$K\subseteq(\sg)^{-1}[^\ast(3/4,1]]$  and  $\sg(\alpha)$  is a positive infinitesimal. 
Hence, $0\le\hg(a) = \text{st}(\sg(\alpha))\le 1/4$  and 
$3/4\le \hg(q(\gamma))=\text{st}(\sg(\gamma))\le 1$  for all  $\gamma\in q^{-1}[F]$. Thus,  
$\hg(a)\notin\text\,\text{cl}\; g[F]$,  as required.   $\tr$

\vskip 0,5 truecm\noindent
{\bf (6.5) Theorem}: $(\hX,\hT)$ coincides with the Hewitt  realcompactification  $\nu X$
of  $(X,T)$ (L. Gillman and M. Jerison [3]).  

\vskip 0,3 truecm\noindent
{\bf Proof}: The space  $(\hX,\hT)$  is realcompact and completely regular, by Theorem (6.2) and 
Theorem (6.4), respectively. Then,  $q[X]$  is a dense subset of   $\hX$  and every  
$f\in C(X,\R)$  has a unique continuous extension  $\hf$  to  $\hX$, by Theorem (3.5) and 
Theorem (3.6), respectively (both applied for  $\Phi = C(X,\R)$). These properties characterize  
$\nu X$.   $\tr$

\vskip 0,5 truecm
E. Hewitt has shown that the real maximal ideals   $M_p$  of  $C(X,\R)$  are uniquely 
determined by points  $p$  in  $\nu X$   by ``$f\in M_p$  iff  $\hf(p) = 0$", where  $\hf$  
denotes the unique extension of  $f$  to  $\nu X$.  

\vskip 0,2 truecm
For completeness we derive this result using the nonstandard methods developed so far:

\vskip 0,5 truecm\noindent
{\bf (6.6) Theorem}:  Let  $M$  be a real  maximal ideal of  $C(X,\R)$. Then, there is a unique 
point   $p$  in  $\nu X$  such that ``$f\in M$  iff  $\hf(p)= 0$".

\vskip 0,3 truecm\noindent
{\bf Proof}: By Lemma (6.1), there is  $\alpha\in\tX$  such that ``$\pi(f) = \ssf(\alpha)$  
for all  $f\in C(X,\R)$"  where  $\pi:C(X,\R)\to C(X,\R)/M = \R$  is the ring homomorphism 
onto  $\R$  determined by  $M$.  Now,  $\hf(q(\alpha)) = \ssf(\alpha)$, i.e. ``$f\in M$  iff  
$\hf(p) = 0$"  for  $p = q(\alpha)$.  The point  $p$  is unique since  $\hX$  is completely 
regular, by Theorem (6.4) and Hausdorff, by Theorem (3.9). The proof is complete.  $\tr$

\vskip 1 truecm
\noindent
{\bf PART III. MONADS AND SEPARATIONS PROPERTIES }

\vskip 0,5 truecm
We study the separation properties of topological spaces such as  $T_0, T_1$, 
{\em regularity, normality, complete regularity, compactness and soberness}
which are characterized in terms 
of monads. Some of the characterizations have already counterparts in the literature on 
nonstandard analysis (but ours are, as a rule, simpler), while others are treated in nonstandard 
terms for the first time. In particular, it seems that the nonstandard characterization of the 
{\em sober spaces} has no counterparts in the nonstandard literature. We also present two new 
characterizations of the compactness in terms of monads similar to but different from A. 
Robinson's famous theorem.

\vskip 0,2 truecm 
We shall use as well the terminology of  (J.L. Kelley [12]) and (L. Gillman and M. Jerison [3]). 
\noindent\\

{\bf 1. Monads and Compactness}

\vskip 0,2 truecm 
A. Robinson proved that a set  $A\subseteq X$ is compact in  $(X,T)$  iff  $\sA$  consists of 
nearstandard points only ([17], Theorem 4.1.13, p. 93). The purpose of this section is to give 
two similar characterizations of compactness in terms of monads which seem to be new in 
the literature on nonstandard analysis.         

\vskip 0,2 truecm
For the definition and the basic properties of the monads the reader should refer to 
(Chapter II, Section 1). As a convenient technique we use the nonstandard compactification  
$(\sX,\,^sT)$  of  $(X,T)$, described in Chapter II, Section 2. Recall that for any 
$H\subseteq X$  we have             
$$
  ^\ast(\text{cl}_X\,H) = \text{cl}_{\sX}\,\sH = \text{cl}_{\sX}\,H
$$
where  $\text{cl}_X$  and  $\text{cl}_{\sX}$  are the closure operators in  $(X,T)$  and  
$(\sX,\,^sT)$,  respectively (Chapter II, Lemma 2.6). Notice that  $\text{cl}_{\sX}$  
coincides with the  $\cF$-monad, in symbols,  $\mu_{\cal F} = \text{cl}_{\sX}$, where  
$\cF$  is the family of all closed sets of  $(X,T)$  (Definition II.1.14). 

\vskip 0,3 truecm\noindent
From Corollary II.1.6, it follows immediately  that
$$
  \bigcup_{\alpha\in\cA}\;\mu(\alpha)\;\subseteq\; \mu(\cA) \leqno{\bf (1.1)}
$$
for any  $\cA\subseteq\sX$.  The next example shows that this inclusion may be proper.

\vskip 0,3 truecm\noindent
{\bf Example}: Let $\N$  be the set of the natural numbers with the discrete topology.   
For any  $n\in\N$  we have  $\mu(n) = \{n\}$,  so that the union of all monads of points 
in  $\N$  is the whole set  $\N$.  On the other hand,  we have   $\mu(\N) = \sN$.

\vskip 0,2 truecm
The next result shows that the equality in (1.1)  holds for subsets  $\cA$  of  $\sX$  
which are compact  in  $(\sX,\,^sT)$.

\vskip 0,3 truecm\noindent
{\bf (1.2) Theorem}:  Let  $(X,T)$  be a topological space and  $(\sX,\,^sT)$  be its 
nonstandard compactification (Definition II.2.4).  If  $\cA\subseteq\sX$  is compact in 
$(\sX,\,^sT)$,  then 
$$
  \bigcup_{\alpha\in\cA}\;\mu(\alpha)\;=\; \mu(\cA) \leqno{\bf (1.3)}
$$

\vskip 0,2 truecm\noindent
{\bf Proof}:  Let  $\alpha\in\mu(\cA)$  and suppose that  $\alpha\notin\mu(\beta)$  
for all  $\beta\in\cA$.  Hence, for any  $\beta\in\cA$  there is  $G_\beta\in T$  such that  
$\beta\in\sG_\beta$  and $\alpha\notin\sG_\beta$. On the other hand, we have, obviously,  
the cover:
$$
  \cA\subseteq\bigcup\{\sG_\beta : \beta\in\cA\}.
$$
Now, by the compactness of  $\cA$,   there exist  $G_{\beta_1},\dots, G_{\beta_n}$  such that  
$$
  \cA\subseteq \{\sG_{\beta_i}:i=1,\dots,n\}\; =\;  
          ^\ast\left(\bigcup\{G_{\beta_i}:i=1,\dots,n\}\right).
$$
By the definition of  $\mu(\cA)$,  we obtain 
$\mu(\cA)\subseteq\,^\ast(\bigcup\{G_{\beta_i} : i = 1,\dots,n\})$,  i.e.  $\alpha\in\sG_{\beta_i}$, 
for some $i$,  which is a contradiction. The proof is complete. $\tr$

\vskip 0,5 truecm\noindent
{\bf (1.4) Corollary}:  Let  $(X,T)$  be a topological space and  $\sX$  be the nonstandard 
extension of  $X$. Then (1.3) holds for any internal subset  $\cA$  of  $\sX$.

\vskip 0,4 truecm\noindent
{\bf Proof}: The internal subsets of $\sX$  are compact in $(\sX,\,^sT)$  (Theorem II. 2.5) 
and the result follows immediately from Theorem (1.2).  $\tr$

\vskip 0,5 truecm\noindent
The next example shows that the equality  (1.3)  may be true for subsets of  $\sX$  which 
are not compact in $(\sX,\,^sT)$.

\vskip 0,5 truecm\noindent
{\bf Example}:  Let  $X = \R$  with the usual topology  $\tau$.  Then 
$$
  \cA = \{n + h : n\in\N,\quad  h\in\sR,\quad h\approx 0\} 
$$
is not compact in  $(\sR,\,^s\tau)$  but (1.3) is still satisfied.

\vskip 0,2 truecm
In contrast with the above, when  $\cA\subseteq X$,  the equality  (1.3)  provides a 
characterization of compactness of  $\cA$  in  $(X,T)$. 

\vskip 0,6 truecm\noindent
{\bf (1.5) Theorem} {\em (Characterization)}: Let  $A\subseteq X$. Then the following  
conditions are equivalent:

  {\bf (i)}   $\;\,A$  is compact in $(X,T)$.

\vskip 0,2 truecm
  {\bf (ii)}  $\,\sA\subseteq \bigcup\limits_{x\in A}\,\mu(x)$.

\vskip 0,2 truecm
  {\bf (iii)} $\bigcup\limits_{x\in A}\,\mu(x) = \bigcup\limits_{\alpha\in\sA}\,\mu(\alpha)$.

\vskip 0,2 truecm
  {\bf (iv)} $\,\bigcup\limits_{x\in A}\,\mu(x) =  \mu(A)$.

\vskip 0,5 truecm\noindent
{\bf Proof}: (i) $\,\Leftrightarrow\,$ (ii) is A. Robinson's theorem mentioned in the 
beginning of this section. 

\vskip 0,3 truecm
  (ii) $\;\Rightarrow\;$ (iii):  Let  $\sA\subseteq\bigcup_{x\in A}\;\mu(x)$.  
So, $\in\sA$  implies  $\alpha\in\mu(x)$  for some  $x\in A$, which implies  
$\mu(\alpha)\subseteq \mu(x)$, by Corollary II.1.6.  
Hence,  $\mu(\alpha)\subseteq\bigcup_{x\in A}\;\mu(x)$  for  all  $\alpha\in\sA$  which implies
$$
  \bigcup\limits_{\alpha\in\sA}\;\mu(\alpha)\quad\subseteq\quad\bigcup\limits_{x\in A}\;\mu(x).
$$
The inverse inclusion is obvious.  

\vskip 0,3 truecm
  (iii)$\;\Rightarrow\;$ (iv):  $\sA$  is an internal subset of  $\sX$  and hence,  
$$
  \bigcup\limits_{\alpha\in\sA}\;\mu(\alpha)\quad = \quad \mu(\sA)\,,
$$ 
by Corollary (1.4) (applied for  $\cA = \sA)$. For the LHS we have
$$
  \bigcup\limits_{\alpha\in\sA}\;\mu(\alpha)\quad = \quad\bigcup\limits_{x\in A}\;\mu(x),
$$ 
by our assumption, and, on the other hand, $\mu(A) =  \mu(\sA)$ (II.1.4), thus, it follows    
$\bigcup_{x\in A}\;\mu(A)\, = \,\mu(A)$,  as required. 

\vskip 0,3 truecm
  (iv)$\;\Rightarrow\;$ (ii):  we have  $\sA\subseteq\mu(A)$,  by the definition of 
$\mu(A)$ (Definition II.1.1),  since  $A\subseteq\sG$  for some  $G\subseteq T$ implies 
$A\subseteq G$, by Corollary I.4.8, which implies  $\sA\subseteq\sG$, by Theorem I.4.2. 
On the other hand,  $\mu(A) = \bigcup_{x\in A}\;\mu(x)$, by assumption, hence, it follows 
$$
  \sA\;\subseteq\;\bigcup\limits_{x\in A}\;\mu(x),
$$ 
as required.   $\tr$\\

\noindent
{\bf 2. Separation Properties and Monads}

\vskip 0,3 truecm
The purpose of the present section is to give characterizations of  separation properties 
like: $T_0$, $T_1$, regularity, normality, complete regularity and soberness in terms of monads. 
Some of the characterizations have counterparts in the literature on nonstandard analysis,  
while others (as the soberness, for example) are treated in nonstandard terms for the first 
time.

\vskip 0,1 truecm 
Two  sets  $A$  and  $B$  will be called  ``comparable" if  ``$A\subseteq B$  or $A\supseteq B$".

\vskip 0,3 truecm\noindent
{\bf (2.1) Theorem}:  Let  $(X,T)$  be a topological space.   Then:

\vskip 0,3 truecm
{\bf (i)}  $\;\,(X,T)$  is a  $T_0$ - space  iff  $x=y\Leftrightarrow \mu(x)=\mu(y)$  for any 
$x,y\in X$.

\vskip 0,2 truecm
{\bf (ii)} $\,(X,T)$  is a  $T_1$ - space  iff  $x=y\Leftrightarrow\mu(x)$  and $\mu(y)$ 
are comparable for any   $x,y\in X$.

\vskip 0,3 truecm\noindent
{\bf Proof}:  (i)  Let  $(X,T)$  be a  $T_0$ - space (J.L. Kelley [12]) and  $x\not= y$. 
Assume that  $x\in G$  but  $y\notin G$  for some   $G\in T$. That is,  $x\in\sG$ and 
$y\notin\sG$  which implies  $y\notin\mu(x)$,  i.e.  $\mu(x)\notin\mu(y)$. The implication 
$(\mu(x)\notin\mu(y))\Rightarrow x\not= y$ is trivial. Assume now, that the condition in 
(i) is valid and let  $x\not= y$.  Without loss of generality, assume  that
$\alpha\in\mu(y) - \mu(x)$.   In other words, there exists  $G\in T$  such that $x\in G$
but  $\alpha \notin\sG$.  Notice now, that  $y\notin G$ (otherwise, $\alpha\in\mu(y)\subseteq\sG$ 
which is a contradiction). Thus,  $(X,T)$  is  $T_0$ .  

\vskip 0,2 truecm
(ii)  Suppose  $(X,T)$  is a  $T_1$ - space (J.L. Kelley [12]) and  $\mu(x)$  and  $\mu(y)$
are comparable, say,  $\mu(x)\subseteq\mu(y)$.  If  $x\not= y$, we have an open set 
$G=X-\{x\}$  with  $x\notin G$  and  $y\in G$. Hence,  $x\notin\mu(y)$  contradicting the 
assumption that  $\mu(x)\subseteq\mu(y)$. Conversely, suppose the comparability property holds. 
If  $(X,T)$  is not  $T_1$, then there exists some  $x$  and  $y$, $y\notin x$,   such  that
$y\in\text{cl}\{x\}$.  Hence  $x\in\mu(y)$,  we have  $\mu(x)\subseteq\mu(y)$  which implies  
$x = y$,   a contradiction. Thus,  $(X,T)$  is a $T_1$ - space.         $\tr$  

\vskip 0,2 truecm
As a consequence of Theorem (2.1), we shall obtain the characterization of  $T_0$-spaces 
given by (A. Robinson [17], Theorem 4.1.9, p.92) and  the characterization of  $T_1$-spaces  
given in (A.E. Hurd and P.A. Loeb [10],  p.114):

\vskip 0,6 truecm\noindent
{\bf (2.2) Theorem}:  The topological space  $(X,T)$  is:

\vskip 0,3 truecm
{\bf (i)} $\;\, T_0$  iff  $x\not= y\Rightarrow \;\text{either} \; x\notin\mu(y)$  or
$y\notin\mu(x)$  for all  $x,y\in X$.

\vskip 0,2 truecm
{\bf (ii)}  $\,T_1$  iff   $x\not=y\Rightarrow x\notin\mu(y)$  and $y\notin\mu(x)$ 
for all  $x,y\in X$.

\vskip 0,5 truecm\noindent
{\bf Proof}:  (i) Let  $(X,T)$  be  $T_0$  and suppose  $x\not= y$.  Then  $\mu(x)\not=\mu(y)$,
 hence, by Corollary (II.1.6),  $x\notin\mu(y)$  or  $y\notin\mu(x)$. Conversely, suppose  
$x\not= y$. Then,  $x\notin\mu(y)$ or  $y\notin\mu(x)$,  hence  $\mu(x)\not=\mu(y)$.

\vskip 0,3 truecm
    (ii)  Let  $(X,T)$  be  $T_1$  and suppose  $x\not= y$. Then  $\mu(x)$  and  $\mu(y)$
 are not comparable. If  $x\in\mu(y)$  or  $y\in\mu(x)$,  then  $\mu(x)\subseteq\mu(y)$
 or   $\mu(y)\subseteq\mu(x)$  which is a contradiction. Conversely, suppose  $x\not= y$. 
Then we have both  $x\notin\mu(y)$  and  $y\notin\mu(x)$  which implies, by Corollary (II.1.6),  
``not  $\mu(x)\subseteq\mu(y)$"  and ``not  $\mu(y)\subseteq\mu(x)$".   $\tr$

\vskip 0,3 truecm
Concerning  $T_2$-spaces, we recall that  a topological space  $(X,T)$  is Hausdorff  
(or a $T_2$ - space) if and only if  $\mu(x)\cap\mu(y) = \emptyset$  for all  $x,y\in X$
$x\not= y$ (A. Robinson [17], Theorem 4.1.8, p. 92). A related notion is that of weakly-Hausdorff: 
$(X,T)$  is called ``weakly Hausdorff"  if  for any   $x\in X$  and  any open neighbourhood  
$G$  of  $x$  the point  $x$  is separated from  all points  $y\in X - G$. The corresponding 
nonstandard characterization is:  $(X,T)$  is  weakly  Hausdorff   if and only if  for  any  
$x,y\in X$  either   $\mu(x) = \mu(y)$,  or  $\mu(x)\cap\mu(y) = \emptyset$  
(K.D. Stroyan and W.A.J. Luxemburg [22],  p.199). 

\vskip 0,2 truecm
We now characterize regularity and normality.

\vskip 0,3 truecm\noindent
{\bf (2.3) Theorem}: Let  $(X,T)$  be a topological space. Then:

\vskip 0,3 truecm
{\bf (i)} $\;\,(X,T)$ is normal iff $F_1\cap F_2=\emptyset\Rightarrow\mu(F_1)\cap\mu(F_2)=\emptyset$
for  any  two closed sets  $F_1, F_2\subseteq X$.

\vskip 0,2 truecm
{\bf (ii)} $\,(X,T)$ is regular iff $\alpha\notin\mu(x)\Rightarrow\mu(\alpha)\cap\mu(x)=\emptyset$
for any  $\alpha\in\sX$  and  $x\in X$.

\vskip 0,3 truecm\noindent
{\bf Proof}: (i) is a version (in terms of closed sets) of  $A$. Robinson's characterization 
of the normality given in ([17], Theorem 4.1.12, p.93) without proof. For completeness we present 
a simple proof: The condition is, obviously,  necessary.  Suppose that  $(X,T)$  is not normal. 
Then, there exist two closed disjoint sets  $F_1$ and $F_2$ such that $U_1\cap U_2\not=\emptyset$  
for all $U_1, U_2\in T$  such that  $F_1\subseteq U_1$  and  $F_2\subseteq U_2$. 
By the saturation principle, we obtain
$$
  \mu(F_1)\cap\mu(F_2) = \bigcap\{\sU_1\cap\sU_2 :
  F_1\subset U_1,F_2 \subset U_2,\, U_1,U_2 \in T\}\not= \emptyset.
$$

  (ii)  Let  $(X,T)$  be a regular space and  $\alpha\in\sX$  and  $x\in X$  be such that 
$\alpha\notin\mu(x)$.  Then,  there is  $G\in T$   such that  $x\in G$  and  $\alpha\notin\sG$. 
By regularity, there is   $U\in T$   such  that  $x\in U$  and  
$\text{cl}_X\,U\subseteq G$.  
We have  $\mu(x)\subset\sU$  and  also  $\mu(\alpha)\subset\,^\ast(X-\text{cl}_X\,U)$  
since  $\alpha\in\sX - \sG = \,^\ast(X-G)\subset\,^\ast(X-\text{cl}_X\,U)$. That is 
$\mu(\alpha)\cap\mu(x) = \emptyset$.  Conversely, suppose that  $(X,T)$  is not regular. 
We show now that there exist  $\alpha\in\sX$  and  $x\in X$  such that $\alpha\notin\mu(x)$
and  $\mu(\alpha)\cap\mu(x)\not=\emptyset$. Indeed, since  $X$  is not regular, there are  
$x\in X$  and   $G\in T$  such that  $x\in G$  and                   
$$
  \text{cl}_X\,H\cap(X-G)\not=\emptyset
$$
for all  $H\in T$  containing  $x$.  By the {\em Saturation Principle} (Chapter I, Section 2, 
Axiom 3), there exists  $\alpha$  such that
$$
  \alpha\in\bigcap\{\,^\ast(\text{cl}_X\,H) : H\in T,\; x\in H\} - \sG.
$$
Since  $\alpha\in\,^\ast(\text{cl}_X\,H) = \text{cl}_{\sX}\,\sH$  then  $\sO\cap\sH\not=\emptyset$  
for all  $O,H\in T$   such that  $\alpha\in\sO$  and  $x\in H$.  Also we have  
$\alpha\notin\mu(x)$,  since  $\alpha\notin\sG$. Using the saturation principle again, we obtain
$$
   \mu(\alpha)\cap\mu(x) = \bigcap\{\sO\cap\sH : O,\;H\in T,\; \alpha\in\sO,\; x\in H\}\not= \emptyset.
$$
The proof is complete.   $\tr$

\vskip 0,5 truecm
As a consequence of Theorem (2.3), we shall obtain the description of regularity given 
in (A. Robinson [17], Theorem 4.1.11, p.93):

\vskip 0,5 truecm\noindent
{\bf (2.4) Theorem}: The topological space  $(X,T)$  is regular iff  
$x\notin F\Rightarrow \mu(x)\cap\mu(F) = \emptyset$  for any  $x\in X$  and any closed  $F\in X$.     

\vskip 0,5 truecm\noindent
{\bf Proof}:  Let the condition hold  and let  $\alpha\in\sX$  and  $x\in X$  and 
$\alpha\notin\mu(x)$. Then there exists  $G$ open  such that  $x\in G$  and  $\alpha\notin\sG$. 
With  $F = X-G$, we have $x\notin F$, so $\mu(\sF)\cap\mu(x) = \emptyset$.  Since $\sF$  is an 
internal set it follows  $\mu(\beta)\cap\mu(x) = \emptyset$  for all $\beta\in\sF$, 
by Corollary (1.4).  Since  $\alpha\in\sF$,  we have  $\mu(\alpha)\cap\mu(x) = \emptyset$.  
Conversely,  let $(X,T)$  be a regular space and   $x\in X$, $F\subseteq X$  be closed  and 
$x\notin F$.  Since  $F$  is closed,  we have  $\sF\cap\mu(x)=\emptyset$, i.e. $\alpha\notin\mu(x)$,
for all  $\alpha\in\sF$ (Theorem, II.1.8). By Theorem (2.3),  $\mu(\alpha)\cap\mu(x)=\emptyset$ 
for all  $\alpha\in\sF$  which implies
$$
  (\bigcup_{x\in F}\,\mu(x))\;\cap\; \mu(x) = \emptyset.
$$
Using (1.4) and (II.1.4),  we obtain  $\mu(F)\cap\mu(x) = \emptyset$, as required.  $\tr$

\vskip 0,5 truecm
Let  $(X,T)$  be a topological space.  Recall that a closed subset  $A\subseteq X$  is  called  
``irreducible"  if  for any  closed subsets  $F_1,\,F_2\subseteq X$  the equality   
$A = F_1\cup F_2$   implies  either  $A = F_1$, or  $A = F_2$.  The space  $(X,T)$  is called  
``sober"  if every irreducible closed subset  $A\subset X$  is  of  the type  $A=\text{cl}_X\{x\}$  
for some  $x\in X$ (see  R.E. Hoffmann [9]  for general reference). 

\vskip 0,5 truecm\noindent
{\bf (2.5) Theorem}:  Let  $A\subset X$  be  a closed set and  $x\in A$.  The following are 
equivalent:

{\bf (i)} $\quad  A = \text{cl}_X \{x\}$.

\vskip 0,2 truecm
{\bf (ii)} $\; \mu(x) = \bigcap\limits_{\alpha\in\sA}\;\mu(\alpha)$.

\vskip 0,2 truecm
{\bf (iii)} $\,\mu(x) = \bigcap\limits_{x\in A}\; \mu(x)$.

\vskip 0,8 truecm\noindent
{\bf Proof}: (i)$\,\Rightarrow\,$(ii): Suppose $A=\text{cl}_X\{x\}$. Since $x\in A$, we get immediately,
$$
  \bigcap_{\alpha\in\sA}\;\mu(\alpha)\subseteq \mu(x)\,.
$$
To show the reverse inclusion, let $\alpha\in´\sA$  and  $G\in T$  be such that 
$\alpha\in\sG$. Then  $A\cap\sG\not=\emptyset$  since  $\sA = \text{cl}_{\sX}\,A$.  
Hence  $\sG\cap\text{cl}_{\sX}\{x\}\not=\emptyset$, so that $x\in G$. Hence $\mu(x)\subseteq\sG$, 
so that $\mu(x)\subseteq\mu(\alpha)$.  This implies
$$
  \mu(x) \subset \bigcap_{\alpha\in\sA}\,\mu(\alpha).
$$

\vskip 0,2 truecm
(ii) $\,\Rightarrow\,$ (iii):  $\mu(x) = \bigcap\limits_{\alpha\in\sA}\,\mu(\alpha)\subseteq
     \bigcap\limits_{x\in A}\,\mu(x)\subseteq\mu(x)$  since  $x\in A$.

\vskip 0,2 truecm
(iii) $\,\Rightarrow\,$ (i): We have  $x\in\mu(\alpha)$  for  all  $a\in A$, so that  
$a\in\text{cl}_X\{x\}$, hence  $A\subseteq\text{cl}_X\{x\}$. On the other hand,  
$\text{cl}_X\{x\}\subseteq A$, since  $A$ is closed and  $x\in A$. The proof is complete.  $\tr$

\vskip 0,6 truecm\noindent
{\bf Example}:  The following example shows that the condition  $x\in A$  for the point  $x$
 is necessary in the above proposition:  $X = \{0, 1, 2\}$,  
$T = \{\emptyset,\{0\}, \{0,1\}, \{0,2\}, X\}$,  $A = \{1,2\},\,x = 0$. Then we have  
$\sX = X$,  $\sT = T,  \sA = A$,  $\mu(0) = \{0\}$,  $\mu(1) = \{0,1\}$, $\mu(2) = \{0,2\}$. 
So we get  $\bigcap_{a\in A}\,\mu(a)=\{0\} = \mu(0)$  in spite of  $x\notin A$. 

\vskip 0,5 truecm\noindent
{\bf (2.6) Definition} {\em (Partial Order in  $\sX$)}:  Let $\alpha, \beta\in\sX$.  Then:

\vskip 0,3 truecm
{\bf (i)} $\;\,\alpha\le \beta$  if  $\mu(\alpha)\subseteq\mu(\beta)$.

\vskip 0,2 truecm
{\bf (ii)} $\,A$  set  $S\subseteq\sX$  is downward directed if for any  $\alpha$, $\beta\in S$
           there is  $\gamma\in S$  such that  $\gamma\le\alpha$  and  $\gamma\le\beta$.

\vskip 0,3 truecm\noindent
{\bf Note}:  If  $x,y\in X$, then  $x\le y$ if and only if  $y\in\text{cl}_X\{x\}$,  
so the order defined above on  $\sX$  is the inverse of the specialization order on  $(X,T)$
(see e.g.  R. E. Hoffmann [9]).    

\vskip 0,3 truecm\noindent
{\bf (2.7) Theorem}: Let $(X,T)$  be a topological space and  $A\subseteq X$  be closed. 
Then  $A$  is irreducible  iff  $\sA$   is a downward directed set.  

\vskip 0,3 truecm\noindent
{\bf Proof}:  Suppose  $\sA$  is downward directed. If  $A$  is not  irreducible,  then  there 
are closed sets  $F_1, F_2\subset A$  such that  $A = F_1 \cup F_2$  and  $A - F_1\not=\emptyset$,
$A - F_2\not=\emptyset$.  Let  $a_i\in A - F_i$.  By assumption, there is  $\gamma\in\sA$  such 
that $\mu(\gamma)\subseteq\mu(a_1)\cap\mu(a_2)$. Now $\gamma\in\mu(a_1)$, so $\gamma\in\sX-\sF_1$;
similarly,  $\gamma\in\sX-\sF_2$. Hence  $\gamma\in(\sX - \sF_1)\cap (\sX - \sF_2) = \sX - \sA$, 
which is a contradiction. Conversely,  suppose  $A$  is an irreducible closed set. Let  
$\alpha, \beta\in\sA$.  We show that  $\sA\cap\mu(\alpha)\cap\mu(\beta)\not= \emptyset$,  
from which the result follows. For  any  $G_i$,  open, such that  $\alpha\in\sG_1$, 
$\beta\in\sG_2$, we have  $A\cap G_1 \cap G_2\not=\emptyset$: otherwise,  $A=(A-G_1)\cup(A-G_2)$ 
so that  $A\subset A-G_1$  or  $A\subseteq A-G_2$. If  $A\subseteq A-G_1$, then  
$G_1\cap A = \emptyset$, so that  $\sG_1\cap\sA = \emptyset$, which contradicts our assumption 
concerning  $\alpha$.  Similarly,  $A\subset A-G_2$  is impossible.  Hence  
$A\cap G_1\cap G_2\not=\emptyset$  for  all  open  $G_i$  such that  $\alpha\in\sG_1$,  
$\beta\in\sG_2$. By  the saturation principle,  
$$
  \sA\cap\mu(\alpha)\cap\mu(\beta)\not= \emptyset
$$   
which finishes the proof.    $\tr$

\vskip 0,5 truecm\noindent
{\bf (2.8) Theorem}:  The topological space $(X,T)$  is {\em sober}  iff  for any closed set 
$A\subseteq X$   such that  $\sA$  is downward directed,  $\sA$  has a smallest element in $A$. 

\vskip 0,5 truecm\noindent
{\bf Proof}:  Suppose  $(X,T)$  is sober.  Let $A$ be a closed set such that  $\sA$  is downward 
directed. Then  $A$  is irreducible hence,  $A = \text{cl}_X\{x\}$  for some  $x\in A$. 
By  Theorem (2.5), 
$$
  \mu(x)\subseteq\bigcup_{\alpha\in\sA}\;\mu(\alpha), \leqno{\bf (2.9)}
$$
i.e.  $x\le\alpha$  for all  $\alpha$  in  $\sA$.  Conversely,  let  $A$  be an irreducible 
closed set in  $X$. Then  $\sA$  is  downward directed, hence  there is  $x\in A$  such that  
(2.9)  holds  which implies  $\mu(x)=\bigcap\,\{\mu(\alpha) : \alpha\in\sA\}$. By Theorem (2.5),  
we get  $A = \text{cl}_X\{x\}$. The proof is complete.   $\tr$

\vskip 0,5 truecm
Recall that  the space  $(X,T)$  is called functionally separated if  for any  $x,y\in X$,   
$x\not= y$,  there exists a continuous function  $f : (X,T)\to(\R,\tau)$,  where  $\tau$  
is the usual topology on  $\R$,  such that  $f(x)=0$  and   $f(y)=1$. We now characterize 
these spaces in terms of special monads.

\vskip 0,2 truecm
Let  $Z = \{f^{-1}[\{0\}]:f\in C(X,\R)\}$  be the family of the zero sets of  continuous real-
valued functions on  $X$  and let  $\mu_Z$  be the corresponding $Z$-monads (Definition (II.1.14)). 
Then we can  give the following characterization: 
\noindent
{\bf (2.10) Theorem}:  If  $(X,T)$  is a topological space,  then:

\vskip 0,3 truecm
{\bf (i)} $\; (X,T)$  is functionally separated iff  
          $x\not= y\Rightarrow\mu_Z(x)\cap\mu_Z(y) = \emptyset$  for  any  $x,y\in X$.

\vskip 0,2 truecm
{\bf (ii)} The topological space  $(X,T)$  is completely regular iff  
           $x\notin F\Rightarrow\mu_Z(x)\cap\mu_Z(F) = \emptyset$   for any  $x\in X$
           and any closed  $F\subseteq X$.

\vskip 0,3 truecm
{\bf (iii)} $\,(X,T)$  is normal iff  
           $F_1\cap F_2 = \emptyset\Rightarrow\mu_Z(F_1)\cap\mu_Z(F_2) = \emptyset$  
           for any closed subsets  $F_1, F_2$  of  $X$.

\vskip 0,5 truecm\noindent
{\bf Proof}:  We shall present the proof of  (ii)  only; the others  are  proved similarly. 
Suppose  $(X,T)$  is completely regular (not necessarily Hausdorff) and let  $x\notin F$. 
Then there is a continuous function  $f: (X,T)\to(I,\tau)$  such that  $f(x)=0$  and  
$F\subseteq f^{-1}[\{1\}]$  where $I = [0,1]$  and  $\tau$  is the usual topology of  $I$.  
Let us set  $Z_0 = [0,1/4]$  and  $Z_1 = [3/4,1]$.  Now we have  
$\mu_Z(x)\subseteq\ssf^{-1}[\sZ_0]$, $\;\mu_Z(F)\subseteq\sg^{-1}[\{0\}]\subseteq\ssf^{-1}[\sZ_1]$ 
for  $g = f-1$ and hence, $\mu_Z(x)\,\cap\,\mu_Z(F) = \emptyset$. Conversely,  suppose that the 
condition holds and let  $x\notin F$. By assumption, $\mu_Z(x)\,\cap\,\mu_Z(F)=\emptyset$. 
Hence, by Lemma (II.1.15), there are zero sets  $Z_0$  and  $Z_1$  in  $X$  such that  
$x\in Z_0$,  $\,F\subseteq Z_1$  and  $Z_0\,\cap\,Z_1 = \emptyset$. Now there exists  
$f:(X,T)\to(I,\tau)$  such that  $Z_0\subseteq f^{-1}[\{0\}]$, $\;Z_1\subseteq f^{-1}[\{1\}]$.  
Then $f(x)=0$  and  $F\subseteq f^{-1}[\{1\}]$, hence $(X,T)$  is completely regular.  
The proof is complete.     $\tr$
\noindent\\

{\bf 3. Topological Applications}

\vskip 0,3 truecm
As an application of the previous characterizations, we present simple proofs of well 
known separation properties for topological spaces.

\vskip 0,2 truecm
Let  $(X,T)$  be a topological space and define an equivalence relation on  $X$  by:  $x\sim y$ 
if  $\text{cl}_X\{x\} = \text{cl}_X\{y\}$.  Let  $q$  be the quotient mapping from  $X$  onto  
$X/_\sim = \tX$, and topologize  $\tX$  by:  $V\subseteq \tX$  is a  $\tT$ - neighbourhood of 
$q(x)$  iff  $q^{-1}[V]$ is a neighbourhood of  $x\in X$.  We also have the special property  
$q^{-1}[q[G]] = G$  for all  $G\in T$, so that  $q$ is an open mapping. The space 
$(\tX,\tT)$  is called  the  $T_0$ - reflection of  $(X,T)$ (see, for example,  H. Herrlich [8]).

\vskip 0,3 truecm\noindent
{\bf (3.1) Theorem}: Let  $(X,T)$  be a topological space. Then, $(X,T)$  is weakly Hausdorff  
iff  $(\tX,\tT)$   is Hausdorff. 

\vskip 0,3 truecm\noindent
{\bf Proof}: Suppose  $\tX$  is  Hausdorff.  To show that  $X$  is weakly Hausdorff,  assume  
that  $\mu(x)\not= \mu(y)$  which, by Corollary (II.1.6),  implies either  $x\notin\mu(y)$  
or  $y\notin\mu(x)$.  Hence $y\notin\text{cl}_X\{x\}$  or  $x\notin\text{cl}_X\{y\}$  which 
implies  $q(x)\not= q(y)$.  Since  $\tX$   is Hausdorff, there are open disjoint sets  
$U$, $V$ in  $\tX$  such that  $q(x)\in U$, $q(y)\in V$.  So that,  $x\in q^{-1}[U]$, 
$y\in q^{-1}[V]$  and  $q^{-1}[U]\cap q^{-1}[V]=\emptyset$. Hence, $\mu(x)\cap\mu(y)=\emptyset$.  
Conversely,  assume  $X$  is weakly Hausdorff.  To show that  $\tX$  is Hausdorff, consider  
$x,y$ such that  $q(x)\not= q(y)$,  i.e.  $\text{cl}_X\{x\}\not= \text{cl}_X\{y\}$.   
Then  $\mu(x)\not= \mu(y)$.  So,  $\mu(x)\cap\mu(y) = \emptyset$  by assumption.  Hence, 
there are open disjoint sets  $G,H$  in  $X$  such that  $x\in G$, $y\in H$. But then  
$q(x)\in q(G)$,  $q(y)\in q(H)$ and $q(G)\cap q(H) = \emptyset$.  Hence,  $X$  is Hausdorff  
since  $q(G)$  and  $q(H)$  are open.   $\tr$

\vskip 0,5 truecm\noindent
{\bf Example}:  Let  $\bf G$  be a topological group. Then the closure of the identity  
$\text{cl}\{e\}$  is a normal subgroup of  $\bf G$ . Then the corresponding factor - group  
$\bf G/\text{cl}\{e\}$  is a Hausdorff topological group, so  $\bf G$  is weakly  Hausdorff. 

\vskip 0,5 truecm\noindent 
{\bf (3.2) Theorem}:  If  $(X,T)$  is  $T_0$  and weakly Hausdorff, then it is  Hausdorff.

\vskip 0,5 truecm\noindent
{\bf Proof}:  Suppose  $x\not= y$. Then  $\mu(x)\not= \mu(y)$  since  $(X,T)$  is  $T_0$ 
which implies  $\mu(x)\cap \mu(y) = \emptyset$,  since  $(X,T)$  is weakly Hausdorff.  
So,  $(X,T)$  is Hausdorff.  $\tr$

\vskip 0,5 truecm\noindent
{\bf (3.3) Theorem}:  If  $(X,T)$  is  $T_0$  and regular,  then  $(X,T)$  is Hausdorff.

\vskip 0,5 truecm\noindent
{\bf Proof}:  Let  $x,y\in X$ and  $x\not= y$. Then  $\mu(x)\not= \mu(y)$,  since  $(X,T)$
is  $T_0$,  i.e. we have either  $x\notin\mu(y)$, or  $y\notin\mu(x)$. On the other hand, 
by regularity, we have  $\mu(\alpha)\cap\mu(y) = \emptyset$  for all $\alpha\in\sX$  such 
that  $\alpha\notin\mu(y)$,  in particular,  $\mu(x)\cap\mu(y) = \emptyset$.  
The proof is complete.    $\tr$

\vskip 0,5 truecm\noindent
{\bf (3.4) Theorem}:  If  $(X,T)$  is compact and Hausdorff,  then  $(X,T)$  is regular.

\vskip 0,5 truecm\noindent
{\bf Proof}: Let  $\alpha\in\sX$  and  $x\in X$  be such that  $\alpha\notin\mu(x)$.  
Now, $\alpha\in\mu(y)$  for some  $y\in X$,  since  $X$  is compact. We have  $x\not= y$, 
by the choice of  $x$  and  $y$,  so,  $\mu(x)\cap\mu(y) = \emptyset$, since  $X$  is Hausdorff.  
On the other hand, $\mu(\alpha)\subseteq\mu(y)$, by Corollary (II.1.6), hence, 
$\mu(x)\cap\mu(\alpha) = \emptyset$.     $\tr$

\vskip 0,5 truecm\noindent
{\bf (3.5) Theorem}:  If  $(X,T)$  is compact and regular, then  $(X,T)$  is normal. 

\vskip 0,3 truecm\noindent
{\bf Proof}:  Let  $F_1$  and  $F_2$  be disjoint closed sets of  $X$. Since  $F_1$  is closed 
and  $F_2\subseteq X - F_1$, we have  $\sF_1\cap\mu(x) = \emptyset$   for any  $x\in F_2$.  
Hence, $\alpha\notin\mu(x)$  for  any  $\alpha\in\sF_1$  and any  $x \in F_2$.  By regularity 
of  $(X,T)$,  we have  $\mu(\alpha)\cap\mu(x) = \emptyset$ ((3.3)) and hence,
$$
  \mu(\alpha)\cap\Bigl(\;\bigcup_{\alpha\in F_2}\;\mu(\alpha)\Bigr) = \emptyset
$$
for any  $\alpha\in\sF_1$.  By compactness of  $(X,T)$  and, hence, of  $F_2$ , we obtain  
$\mu(\alpha)\cap\mu(F_2) = \emptyset$  for any  $\alpha\in\sF_1$ , by Theorem (1.5),  
which immediately implies
$$
  \Bigl(\;\bigcup_{\alpha\in F_1}\;\mu(\alpha)\Bigr)\;\cap\;\mu(F_2) = \emptyset.
$$
Since  $F_1$ is also compact, we get  $\mu(F_1)\cap\mu(F_2) = \emptyset$.  
The proof is complete.   $\tr$

\vskip 0,3 truecm\noindent
{\bf (3.6) Theorem}:  If  $(X,T)$  is Hausdorff, then it is sober.

\vskip 0,3 truecm\noindent
{\bf Proof}:  Let  $A$  be a closed set of  $X$  such that  $\sA$  is downward directed.  
Then we have  $A=\{x\}$ for some $x\in A$. For suppose not, i.e. there are  $x,y\in A, x\not= y$,
we obtain  $\mu(x)\cap\mu(y) = \emptyset$,  so $\sA$  cannot be downward directed. 
Hence,  $x$  is the smallest element of  $\sA$  in  $A$.   $\tr$
\noindent\\

{\bf 4.  Separation Properties of   $(^*X, {^sT})$}

\vskip 0,3 truecm
In this section we apply the results  established so far to study the separation properties 
of space  $(\sX,\,^sT)$ ( Chapter II, Section 2). Our interest in the space  $(\sX,\,^sT)$ 
arises from the importance of this space for compactifications and completions of topological spaces, 
demonstrated in Chapter II of this text, as well as its importance for compactifications  of 
ordered topological spaces (S. Salbany and T. Todorov [19]-[20]). 

\vskip 0,3 truecm\noindent
{\bf (4.1) Theorem}:  $(\sX,\,^sT)$  is normal  iff  $(X,T)$  is normal.

\vskip 0,3 truecm\noindent
{\bf Proof}:  Assume that  $(\sX,\,^sT)$  is normal. Let  $F_1$,  $F_2$   be disjoint closed 
sets of  $(X,T)$.  Then  $\sF_1$  and  $\sF_2$  are disjoint closed sets of  $(\sX,\,^sT)$  so, 
by assumption, they can be included in disjoint open sets with disjoint closures. Restricting 
such open sets to  $X$  provides two disjoint open sets  $G_1$, $G_2$  in  $(X,T)$  whose 
closures in  $(X,T)$  are disjoint  and  $F_i\subseteq G_i$ . Conversely, let  $(X,T)$   
be normal  and  let  $A,B\subseteq\sX$  be disjoint closed subsets of  $(\sX,\,^sT)$.  
Now  $A = \mu_{\cF}(A)$  and  $B = \mu_{\cF}(B)$, since  $\mu_{\cF} = \text{cl}_{\sX}$  so, 
by Lemma (II.1.15),  $A\subseteq\sF_1$  and   $B\subseteq\sF_2$, $\;\sF_1\cap\sF_2=\emptyset$, 
for some disjoint closed sets  $F_1$  and  $F_2$  in  $(X,T)$.  But then, by assumption, 
there are open sets  $G_1$  and  $G_2$  of  $(X,T)$  such that 
$F_1\subseteq G_1\subseteq X-G_2\subseteq X-F_2$.  Hence,  
$\sF_1\subseteq\sG_1\subseteq\sX-\sG_2\subseteq\sX-\sF_2$,  so that  $(\sX,\,^sT)$  is normal.  
The proof is complete.  $\tr$

\vskip 0,3 truecm\noindent
{\bf (4.2) Theorem}:  $(\sX,\,^sT)$  is regular  iff  every open set in $(X,T)$  is closed. 

\vskip 0,3 truecm\noindent
{\bf Proof}: Suppose  $(X,T)$  has the stated property and  $\alpha\in\sG$  for some open set 
$G$ in  $(X,T)$.  Since  $G$  is open and closed,  we have  $\sG$  is open and closed  in   
$(\sX,\,\hsT)$, so  $(\sX,\,^sT)$  is regular.  Conversely, suppose  $(\sX,\,^sT)$  is regular.  
Let  $G$  be an open subset of  $X$. Suppose  $G$  is not closed, so  there is  
$x\in\text{cl}_X\,G-G$. For each open neighbourhood $H$ of $x$, we have $G\cap H\not= \emptyset$,  
so the family $\{\sG\cap\sH : H\in T,\; x\in H\}$  has the finite intersection property.  
By the saturation principle, there is a point  $\alpha$  such that
$$
  \alpha\in\bigcap\{\sG\cap\sH : H\in T,\; x\in H\} = 
   \left(\bigcap\{\sH : H\in T,\; x\in H\}\right)\cap\sG.
$$
By regularity, there is $U$, open in $(X,T)$, such that  
$\alpha\in\sU\subseteq\text{cl}_{\sX}\,\sU\subseteq\sG$.  But then, 
$$
  x\in(\sX - \text{cl}_{\sX}\sU)\cap X
$$
since  $x\notin\sG$,  hence  $x\in W = X - \text{cl}_X\,U$.  Thus  $\sU\cap\sW = \emptyset$'  
(as $U\cap W = \emptyset$),  which contradicts  $\alpha \in\sW$  whenever   $W\in T$  and 
$x\in W$. The proof is complete.     $\tr$

As a consequence of the above we have:

\vskip 0,3 truecm\noindent
{\bf (4.3) Theorem}:  Let  $D$  be the discrete topology on  $\N$. Then  $(\sN,\,^sD)$  
is not a  $T_0$ - space.

\vskip 0,3 truecm\noindent
{\bf Proof}:  If  $(\sN,\,^sD)$  were  $T_0$,  then it would be  $T_2$, since  $(\sN,\,^sD)$ 
is regular. Then, since every bounded continuous real valued function on  $(\N,D)$  admits 
a continuous extension to  $(\sN,\,^sD)$  and  $(\N,D)$  is dense in  $(\sN,\,^sD)$  
(Theorem II.2.5), it follows that  $(\sN,\,^sD)$  is  the Stone - \v Cech  compactification  
$\beta(\N,D) = \beta\N$  of  $(\N,D)$.  It is well known that this is impossible 
(see  A. Robinson [18], p. 582)  or  (K.D. Stroyan and W.A.J. Luxemburg [22], 
(8.1.6), (8.1.7). (9.1)).    $\tr$

\vskip 0,3 truecm\noindent
{\bf (4.4) Corollary}: There is no topology  $T$  on $\N$  for which  $(\sN,\,^sT)$  is  a $T_0$  - space.
Proof:  Suppose the contrary, i.e. that  $(\sN,\,^sT)$  is $T_0$  for some topology  $T$  on  $\N$. 
The identity map   $i:(\N,D)\to(\N,T)$   is continuous,   hence so is  its nonstandard extension:                       
$$
  ^\ast i:(\sN,\,^sD)\to(\sN,\,^sT)
$$  
(Theorem II.2.7). Since  $^\ast i$  is injective  and  $(\sN,\,^sT)$  is  $T_0$,  it follows that  
$(\sN,\,^sD)$  is  $T_0$, which is a contradiction. The proof is complete.        $\tr$

\vskip 0,3 truecm\noindent
{\bf (4.5) Theorem}:  $(\sX,\,^sT)$  is a  $T_0$ - space iff  $X$  is finite.

\vskip 0,3 truecm\noindent
{\bf Proof}: If  $X$  is infinite and  $(\sX,\,^sT)$  is a  $T_0$ - space,  then  $X$  has a 
countable subset  $\N\subseteq X$  with relative topology,  also denoted by   $\,^sT$,  
such that  $(\sN,\,^sT)\subseteq (\sX,\,^sT)$.  Thus  $(\sN,\,^sT)$  is  a  $T_0$ - space, 
which is impossible.  The converse is clear.      $\tr$

\vskip 0,2 truecm
This startling result should not be regarded as indicating that  $(\sX,\,^sT)$  is only 
interesting when  $X$  is finite but rather that topologies and sets of points should be 
considered as they occur in Nature. In particular,  $\sN$  has many points which allow 
the extension of natural numbers and their operations but few open sets in the standard 
topology.  However, there are  still  enough open sets to obtain  $\beta\N$  as a quotient 
space of  $(\sN,\,^sT)$ (Chapter II, Section 4). 

\end{document}